\documentclass[12pt]{article}

\usepackage{amssymb}
\usepackage{amsmath}
\usepackage[all]{xypic}
\usepackage[frenchb]{babel}
\usepackage{a4wide}

\title{Th\'eor\`eme de Dobrowolski-Laurent pour les extensions ab\'eliennes sur une courbe elliptique \`a multiplication complexe}

\author{Nicolas Ratazzi}
\date{}

\begin{document}

\maketitle

\renewcommand{\thefootnote}{}
\footnote{\textit{Email address : }ratazzi@math.jussieu.fr}
\renewcommand{\thefootnote}{1}

\newcounter{ndefinition}[section]
\newcommand{\defi}{\addtocounter{ndefinition}{1}{\noindent \textbf{D{\'e}finition \thesection.\thendefinition.\ }}}
\newtheorem{rem}{Remarque }[section]
\newtheorem{lemme}{Lemme}[section]
\newtheorem{conj}{Conjecture}[section]
\newcounter{nex}[section]
\newcommand{\ex}{\addtocounter{nex}{1}{\noindent \textit{Exemple} \thesection.\thenex\ }}
\newtheorem{prop}{Proposition} [section]
\newtheorem{cor}{Corollaire} [section]
\newtheorem{theo}{Th{\'e}or{\`e}me} [section]
\newcommand{\demo}{\noindent \textit{D{\'e}monstration} : }

\noindent\hrulefill

\noindent \textbf{Abstract :} Let $E/K$ be an elliptic curve with complex multiplication and let $K^{\textnormal{ab}}$ be the Abelian closure of $K$. We prove in this article that there exists a constant $c(E/K)$ such that : for all point $P\in E(\overline{K})\backslash E_{\textnormal{tors}}$, we have 
\[\widehat{h}(P)\geq\frac{c(E/K)}{D}\left(\frac{\log \log 5D}{\log 2D}\right)^{13},\]
\noindent where $D=[K^{\textnormal{ab}}(P):K^{\textnormal{ab}}]$. This result extends to the case of  elliptic curves with complex multiplication the previous result of Amoroso-Zannier \cite{AZ} on the analogous problem on the multiplicative group $\mathbb{G}_m$, and generalizes to the case of extensions of degree D the result of Baker \cite{baker} on the lower bound of the N\'eron-Tate height of the points defined over an Abelian extension of an elliptic curve with complex multiplication. This result also enables us to simplify the proof of a theorem of Viada \cite{viada}.

\vspace{.3cm}



\noindent\hrulefill

\section{Introduction}
\noindent Soit $K$ un corps de nombres. En notant $\widehat{h}$ la hauteur de N\'eron-Tate sur une courbe elliptique $E/K$ et en notant $K^{\textnormal{ab}}$ la cl\^oture ab\'elienne de $K$, on montre dans cet article les deux r\'esultats suivants :

\vspace{.3cm}

\begin{theo}\label{thmoi} Si $E/K$ est une courbe elliptique \`a multiplication complexe, il existe une constante $c(E/K)$ strictement positive, telle que 
\[ \forall P\in E(\overline{K})\backslash E_{\textnormal{tors}},\ \ \widehat{h}(P)\geq \frac{c(E/K)}{D}\left(\frac{\log \log 5D}{\log 2D}\right)^{13},\]
\noindent o\`u $D=[K^{\textnormal{ab}}(P):K^{\textnormal{ab}}]$.
\end{theo}

\vspace{.3cm}

\begin{theo}\label{th2}Soient $c_0>0$ et $E/K$ une courbe elliptique \`a multiplication complexe. Il existe une constante strictement positive $c(E/K,c_0),$ telle que : pour toute extension ab\'elienne $F/K$ et pour tout point $P\in E(\overline{K})\backslash E_{\textnormal{tors}}$ v\'erifiant $D=[F(P):F]$, si le nombre de nombres premiers qui se ramifient dans $F$ est born\'e par $c_0\left(\frac{\log 2D}{\log\log 5D}\right)^2$, alors on a l'in\'egalit\'e 
\[ \widehat{h}(P)\geq \frac{c(E/K,c_0)}{D}\left(\frac{\log\log 5D}{\log 2D}\right)^3.\]
\end{theo}

\vspace{.3cm}

\noindent On voit qu'en imposant une contrainte sur l'\'etendue de la ramification dans l'extension ab\'elienne (th\'eor\`eme \ref{th2}), on obtient une g\'en\'eralisation d'un pr\'ec\'edent r\'esultat de Laurent \cite{laurent} (cf. le th\'eor\`eme \ref{laurent} plus loin). Dans le cas g\'en\'eral (th\'eor\`eme \ref{thmoi}), sans imposer aucune condition, on obtient une minoration optimale aux puissances de log pr\`es, avec un exposant l\'egerement d\'egrad\'e par rapport au cas classique : on a comme puissance de log un exposant $13$ au lieu d'un exposant $3$ ; toutefois cet exposant $13$ est le m\^eme que dans le cas multiplicatif d\^u \`a Amoroso-Zannier \cite{AZ} (cf. th\'eor\`eme \ref{thaz} plus loin). Ce th\'eor\`eme \ref{thmoi}, dans le cas des courbes elliptiques \`a multiplication complexe, g\'en\'eralise au cas $D$ quelconque un pr\'ec\'edent r\'esultat de Baker \cite{baker} (cf. th\'eor\`eme \ref{bakersilver} plus loin). Nous donnons \`a la fin de l'introduction une application de notre th\'eor\`eme \ref{thmoi}.

\vspace{.3cm}

\noindent Ce type de probl\`eme remonte aux travaux de Lehmer dans les ann\'ees 1930 : soit $x\in \mathbb{G}_m(\overline{\mathbb{Q}})\backslash\mu_{\infty}$ un nombre alg\'ebrique qui n'est pas une racine de l'unit\'e. On sait par un th\'eor\`eme de Kronecker que sa hauteur logarithmique absolue $h(x)$ est strictement positive. En 1933 Lehmer \'enonce la c\'el\`ebre conjecture

\vspace{.3cm}

\begin{conj}\label{lehmer}\textnormal{\textbf{(Probl\`eme de Lehmer)}} Il existe une constante $c>0$ telle que 
\[\forall x\in  \mathbb{G}_m(\overline{\mathbb{Q}})\backslash\mu_{\infty},\ \ \ h(x)\geq\frac{c}{D},\]
\noindent o\`u $D=[\mathbb{Q}(x):\mathbb{Q}]$.
\end{conj}

\vspace{.3cm}

\noindent Plus exactement, Lehmer se pose plut\^ot la question inverse : est-il possible de contredire cet \'enonc\'e ?

\vspace{.3cm}

\noindent C'est en 1979 , avec le th\'eor\`eme de Dobrowolski \cite{dob}, qu'est obtenu un r\'esultat optimal \`a des puissances de log pr\`es, en direction de cette conjecture :

\vspace{.3cm}
 
\begin{theo}\label{dobro}\textnormal{\textbf{(Dobrowolski)}} Il existe une constante $c>0$ telle que 
\[\forall x\in  \mathbb{G}_m(\overline{\mathbb{Q}})\backslash\mu_{\infty},\ \ \ h(x)\geq\frac{c}{D}\left(\frac{\log\log 3D}{\log 2D}\right)^3,\]
\noindent o\`u $D=[\mathbb{Q}(x):\mathbb{Q}]$.
\end{theo}

\vspace{.3cm}

\noindent Peu de temps apr\`es, Laurent a \'etendu, dans son article \cite{laurent}, la conjecture de Lehmer aux courbes elliptiques sur un corps de nombres, en rempla\c{c}ant la hauteur sur $\mathbb{G}_m$ par la hauteur de N\'eron-Tate et il a \'etendu le r\'esultat de Dobrowolski au cas des courbes elliptiques $E/K$ \`a multiplication complexe. 

\vspace{.3cm}

\begin{theo}\label{laurent}\textnormal{\textbf{(Laurent)}} Soit $E/K$ une courbe elliptique \`a multiplication complexe. Il existe une constante strictement positive $c(E/K)$ telle que 
\[\forall P\in  E(\overline{K})\backslash E_{\textnormal{tors}},\ \ \ \widehat{h}(P)\geq\frac{c(E/K)}{D}\left(\frac{\log\log 3D}{\log 2D}\right)^3,\]
\noindent o\`u $D=[K(P):K]$.
\end{theo}

\vspace{.3cm}

\noindent Dans les articles \cite{amodvor} et \cite{AZ}, Amoroso-Dvornicich et Amoroso-Zannier ont \'etendu le probl\`eme de Lehmer sur $\mathbb{G}_m$ au cas des extensions ab\'eliennes relatives. Pr\'ecis\'ement, ils \'enoncent la conjecture et d\'emontrent le th\'eor\`eme suivant :

\vspace{.3cm}

\begin{conj}\label{az} \textnormal{\textbf{(Amoroso-Zannier)}} Soit $K$ un corps de nombres. Il existe une cons\-tan\-te strictement positive $c(K)$, telle que 
\[\forall x\in  \mathbb{G}_m(\overline{K})\backslash\mu_{\infty},\ \ \ h(x)\geq\frac{c(K)}{D},\]
\noindent o\`u $D=[K^{\textnormal{ab}}(x):K^{\textnormal{ab}}]$.
\end{conj}

\vspace{.3cm}

\begin{theo}\label{thaz} \textnormal{\textbf{(Amoroso-Zannier)}} Soit $K$ un corps de nombres. Il existe une cons\-tan\-te $c(K)$ strictement positive, telle que 
\[\forall x\in  \mathbb{G}_m(\overline{K})\backslash\mu_{\infty},\ \ \ h(x)\geq\frac{c(K)}{D}\left(\frac{\log\log 5D}{\log 2D}\right)^{13},\]
\noindent o\`u $D=[K^{\textnormal{ab}}(x):K^{\textnormal{ab}}]$.
\end{theo}

\vspace{.3cm}

\noindent Ce th\'eor\`eme \'etend le r\'esultat de Amoroso-Dvornicich qui traitait le cas o\`u $x$ appartenait \`a une extension ab\'elienne de $K$, \textit{i.e.}, le cas $D=1$. C'est pr\'ecis\'ement ce th\'eor\`eme, dans le cas $D=1$, qui a \'et\'e \'etendu aux courbes elliptiques \`a multiplication complexe, ou ayant un $j$-invariant non-entier, par Baker dans \cite{baker}, puis par Silverman \cite{silver} dans le cas des courbes elliptiques sans multiplication complexe. Ainsi pour les courbes elliptiques, on a

\vspace{.3cm}

\begin{theo}\label{bakersilver}\textnormal{\textbf{(Baker-Silverman)}} Soit $E/K$ une courbe elliptique. Il existe une cons\-tante strictement positive $c(E/K)$ telle que 
\[\forall P\in E(K^{\textnormal{ab}})\backslash E_{\textnormal{tors}},\ \ \widehat{h}(P)\geq c(E/K).\]
\end{theo}

\vspace{.3cm}

\noindent L'objectif du pr\'esent article est d'\'etendre le r\'esultat d'Amoroso-Zannier au cas des courbes elliptiques \`a multiplication complexe, g\'en\'eralisant ainsi le r\'esultat de Baker au cas $D$ quelconque. Notons que le th\'eor\`eme \ref{thmoi} r\'epond \`a une conjecture de David dans le cas des courbes elliptiques \`a multiplication complexe : 

\vspace{.3cm}

\begin{conj}\label{david}\textnormal{\textbf{(David)}} Soient $A/K$ une vari\'et\'e ab\'elienne sur un corps de nombres et $\mathcal{L}$ un fibr\'e en droites ample et sym\'etrique sur $A$. Pour tout $\varepsilon>0$, il existe une constante strictement positive $c(A/K,\mathcal{L})$ telle que pour tout point $P\in A(\overline{K})$ qui n'est pas de $\textnormal{End}(A(\overline{K}))$-torsion, on a
\[ \ \ \widehat{h}_{\mathcal{L}}(P)\geq\frac{c(A/K,\mathcal{L})}{D_{\textnormal{tors}}^{\frac{1}{g}+\varepsilon}},\]
\noindent o\`u $D_{\textnormal{tors}}=[K(A_{\textnormal{tors}},P):K(A_{\textnormal{tors}})]$.
\end{conj}

\vspace{.3cm}

\noindent En effet, le th\'eor\`eme \ref{thmoi} \'etant vrai pour $D=[K^{\textnormal{ab}}(P):K^{\textnormal{ab}}]$, il l'est en particulier pour $D=[F(P):F]$ pour toute extension ab\'elienne $F/K$. De plus, quitte \`a remplacer $K$ par une extension de degr\'e born\'e en fonction de $E$, le r\'esultat reste toujours vrai (on ne change que la constante $c(E/K)$). L'extension $H(E_{\textnormal{tors}})/H$ est ab\'elienne pour $H=K(j)$ corps de classes de Hilbert de $E$, ce qui conclut.

\vspace{.3cm}

\noindent Le th\'eor\`eme \ref{thmoi} rend naturel de g\'en\'eraliser la conjecture \ref{az} aux courbes elliptiques :

\begin{conj}Soit $E/K$ une courbe elliptique. Il existe une constante strictement positive $c(E/K)$, telle que 
\[ \forall P\in E(\overline{K})\backslash E_{\textnormal{tors}},\ \ \widehat{h}(P)\geq \frac{c(E/K)}{D},\]
\noindent o\`u $D=[K^{\textnormal{ab}}(P):K^{\textnormal{ab}}]$.
\end{conj}

\vspace{.3cm}

\noindent Le th\'eor\`eme \ref{thmoi} est une premi\`ere \'etape en direction de cette conjecture \ref{david}, au moins dans le cas de multiplication complexe. On peut indiquer bri\`evement un des int\'er\`ets d'un tel r\'esultat. Pour expliquer cela, on introduit quelques notations : on dit qu'une courbe (int\`egre) sur une vari\'et\'e ab\'elienne $A$ est \textit{transverse} si elle n'est contenue dans aucun translat\'e de sous-vari\'et\'e ab\'elienne de $A$ diff\'erente de $A$. Si $X$ est un sous-sch\'ema ferm\'e int\`egre de $A$ et $r$ un entier, alors $Z_{x,0}^{(r)}\subset X(\overline{K})$ est l'ensemble des points pour lesquels il existe un sous-sch\'ema en groupes $G$ de $A$ avec 
\[\textnormal{dim}_PX\cap G\geq\max\left\{1,r-\textnormal{codim }G\right\}.\]
\noindent On dit qu'une vari\'et\'e ab\'elienne simple $A$ est de \textit{type $(g,\delta)$}, si elle est de dimension $g$ et si le rang de $\textnormal{End}(A)=2g/\delta$. Enfin, on note 
\[A^{[r]}=\bigcup_{\textnormal{codim }G\geq r}G(\overline{K}),\]
\noindent o\`u $G$ est un sous-sch\'ema en groupe de codimension indiqu\'ee. Dans son article \cite{remond}, R\'emond prouve :

\vspace{.3cm}

\begin{theo}\label{rem}\textnormal{\textbf{(R\'emond)}} Soit $A$ une vari\'et\'e ab\'elienne sur $\overline{K}$. Nous choisissons une isog\'enie entre $A$ et un produit $A_1^{n_1}\times\cdots\times A_m^{n_m}$ o\`u $m$ est un entier naturel et pour chaque indice $i$ avec $1\leq i\leq m$ la vari\'et\'e ab\'elienne $A_i$ est simple de type $(g_i,\delta_i)$ et $n_i\in \mathbb{N}^*$. Soient $X$ un sous-sch\'ema ferm\'e int\`egre de $A$ et $r,r'$ deux entiers tels que $0\leq r\leq r'\leq\textnormal{dim }A$. Nous supposons que l'une des conditions suivantes est v\'erifi\'ee.

\noindent $(C_1)$ La conjecture (\ref{david}) est vraie.

\noindent $(C_2)$ La vari\'et\'e ab\'elienne $A$ est \`a multiplication complexe et $r'>(1+\sum_{i=1}^m g_i)(r-1)$.

\noindent $(C_3)$ On a l'in\'egalit\'e 
\[r'> \sum_{i=1}^m g_i(n_i+\delta_i)\frac{r-1}{r}.\]
\noindent Alors, pour toute hauteur $h$ associ\'ee \`a un fibr\'e ample $\mathcal{L}$ sur $A$ et tout r\'eel $H$, l'ensemble
\[\left\{P\in \left(X(\overline{K})\backslash Z_{X,0}^{(r)}\right)\cap A^{[r']}\ \ / \ h(P)\leq H\right\}\]
\noindent est fini. Si de plus $X$ est une courbe transverse et $r\geq 2$, alors $X(\overline{K})\cap A^{[r']}$ est fini.
\end{theo}

\vspace{.3cm}

\noindent Notons que notre th\'eor\`eme \ref{thmoi} permet d\'ej\`a de simplifier la preuve du theorem 2. de Viada \cite{viada} suivant :

\vspace{.3cm}

\begin{theo}\textnormal{\textbf{(Viada)}} Soient $E/K$ une courbe elliptique \`a multiplication complexe, $n$ un entier non nul et $C/K$ une courbe transverse dans $E^n$. Pour $r\geq 0$ on consid\`ere les ensembles 
\[S_r(C):=\bigcup_{\textnormal{codim\,}G\geq r}G\cap C(\overline{K})\]
o\`u l'union porte sur les sous-groupes alg\'ebriques $G$ de $E^n$ de codimension au moins $r$. Alors l'ensemble $S_2(C)$ est fini.
\end{theo}

\vspace{.3cm}

\noindent La preuve de Viada est calqu\'ee sur celle de Bombieri, Masser et Zannier \cite{BMZ} dans le cas de $\mathbb{G}_m^n$. Elle utilise le fait que la hauteur des points de $S_1(C)$  est born\'ee. Il s'agit du Theorem 1. du m\^eme article de Viada qui r\'esulte simplement des propri\'et\'es fonctorielles des hauteurs et du th\'eor\`eme du cube pour les vari\'et\'es ab\'eliennes. Ceci \'etant acquis on constate, en appliquant le th\'eor\`eme de Northcott, qu'il suffit alors de montrer que le degr\'e des points de $S_2(C)$ est born\'e. C'est la partie difficile de la preuve. Viada montre ceci en deux \'etapes : la premi\`ere consiste \`a montrer la finitude de l'ensemble $S_3(C)$. La seconde \'etape consiste \`a montrer la finitude de $S_2(C)$ en utilisant un subtil argument cohomologique. Nous montrons ici comment \'eviter cet argument cohomologique en appliquant notre th\'eor\`eme \ref{thmoi}. En fait l'utilisation de ce th\'eor\`eme \ref{thmoi} permet de ramener la seconde \'etape \`a la premi\`ere. Nous expliquons ceci dans la derni\`ere partie de cet article.

\vspace{.3cm}

\noindent Dans la suite (derni\`ere partie except\'ee) on s'attache \`a prouver le th\'eor\`eme \ref{thmoi}. On explique \`a la fin comment le th\'eor\`eme \ref{th2} s'obtient de la m\^eme fa\c{c}on. La preuve est une preuve classique de transcendance \`a deux exceptions pr\`es : on utilise un lemme de Siegel absolu et il y a en fait deux extrapolations selon que l'on est dans une situation avec beaucoup de ramification ou non. Ceci \'etant dit, dans le cas non-ramifi\'e, la preuve suit le sch\'ema initi\'e par Dobrowolski, \`a savoir une extrapolation sur les transform\'es par le morphisme de Frobenius. Dans le cas ramifi\'e, on suit la preuve du cas multiplicatif de \cite{AZ} en utilisant encore des transform\'es par Frobenius. On utilise l'astuce de Laurent \cite{laurent} consistant \`a d\'edoubler les variables pour permettre une plus grande libert\'e dans le choix des param\`etres auxiliaires. La partie \ref{rappel} consiste en des rappels sur la hauteur de N\'eron-Tate et sur les propri\'et\'es dont nous aurons besoin concernant les courbes elliptiques \`a multiplication complexe. La partie \ref{reduc} consiste en une s\'erie de r\'eductions en vue de prouver les th\'eor\`emes \ref{thmoi} et \ref{th2}. La preuve proprement dite se trouve dans les parties \ref{sieg}, \ref{extra} et \ref{conclusion}. 

\vspace{.3cm}

\noindent Dans la preuve on se ram\`ene \`a travailler avec une extension ab\'elienne $F/K$ finie et avec $D=[F(P):F].$ Notons que l'hypoth\`ese ``$F/K$ est ab\'elienne'' sert de mani\`ere cruciale dans les deux \'etapes d'extrapolation : dans l'\'etape o\`u il y a peu de premiers ayant un grand indice de ramification dans $F$, \textit{i.e.} l'\'etape ``quasi-classique'', l'extrapolation se fait gr\^ace au lemme \ref{nonramifie} qui utilise de mani\`ere fondamentale l'hypoth\`ese d'ab\'elianit\'e. Dans l'autre extrapolation, \textit{i.e.} le cas comp\'ementaire o\`u beaucoup de premiers ont un grand indice de ramification dans $F$, l'hypoth\`ese sert \`a fabriquer le groupe $H_p$ du lemme \ref{ramifie} : on utilise pour cela le th\'eor\`eme de Kronecker-Weber.

\section{Hauteur et multiplication complexe}\label{rappel}

\subsection{Hauteur}

\noindent Soient $K$ un corps de nombres de degr{\'e} d, $M_K$ l'ensemble des valeurs absolues (deux \`a deux non {\'e}quivalentes) sur $K$, $M_K^0$ les valeurs absolues ultram\'etriques de $M_K$ normalis{\'e}es par $\vert p\vert_v=p^{-1}$ pour toute place finie $v$ au-dessus du nombre premier $p$ et $M_K^{\infty}$ les valeurs absolues archim\'ediennes de $M_K$. On note $d_v=[K_v:\mathbb{Q}_p]$ le degr{\'e} local et on d{\'e}finit la \textit{hauteur (logarithmique absolue)} sur $\mathbb{P}^n(\overline{\mathbb{Q}})$ par
\[h(x_0:\ldots:x_n)=\frac{1}{d}\sum_{v\in M_K} d_v \log\max_{0\leq
  i\leq n}\vert x_i\vert_v.\]
\noindent Dans cette d{\'e}finition, la renormalisation par $\frac{1}{d}$ sert juste {\`a} faire en sorte que $h(x)$ soit
ind{\'e}pendante du choix du corps $K$ contenant $x$. De plus par la formule du produit, la hauteur est aussi ind{\'e}pendante du choix d'un syst{\`e}me de coordonn{\'e}es projectives.

\vspace{.3cm}

\noindent En plongeant $\mathbb{G}_m^n$ dans $\mathbb{P}^n$ par $(x_1,\ldots,x_n)\mapsto (1:x_1,\ldots:x_n)$ ceci d\'efini \'egalement la hauteur sur $\mathbb{G}_m^n$.

\vspace{.3cm}

\noindent Dans la suite on utilisera \'egalement la hauteur $h_2$ d\'efinie sur $\mathbb{P}^n(\overline{\mathbb{Q}})$ par
\[h_2(x_0:\ldots:x_n)=\frac{1}{d}\left(\sum_{v\in M_K^0} d_v \log\max_{0\leq
  i\leq n}\vert x_i\vert_v+ \sum_{v\in M_K^{\infty}} d_v \log \sqrt{\sum_{0\leq
  i\leq n}\vert x_i\vert^2_v}\right).\] 

\noindent Soit $N$ un entier. On d\'efinit comme le fait Schmidt (voir \cite{schmidt}) la hauteur $h_2$ d'un sous-$\overline{\mathbb{Q}}$-espace vectoriel $S$ alg\'ebrique de dimension $d$ de $\overline{\mathbb{Q}}^{N+1}$ par :
\[h_2(S)=h_2(\mathbf{x}_1\wedge\ldots\wedge \mathbf{x}_d),\]
\noindent o\`u $\mathbf{x}_1,\ldots, \mathbf{x}_d$ est une base de $S$ sur un corps de nombres quelconque sur lequel $S$ est d\'efini.

\subsection{Hauteur de N\'eron-Tate}

\defi Si $E/K$ est une courbe elliptique donn{\'e}e par une {\'e}quation de Weierstrass, on d{\'e}finit la \textit{hauteur} $h\  :\ E(\overline{K})\rightarrow \mathbb{R}^+$ par $h(P):=h(x(P):1)$, o\`u $h(x:y)$ est la hauteur logarithmique absolue sur $\mathbb{P}^1(\overline{K})$ d\'efinie pr\'ec\'edemment.

\vspace{.3cm}

\noindent Cette hauteur v\'erifie un certain nombre de propri\'et\'es. Nous indiquons les plus essentielles, qui nous serviront dans la suite. On renvoie par exemple au livre \cite{hindrysil} Part B pour tout ce qui concerne les hauteurs.

\vspace{.3cm}

\begin{prop} Sur une courbe elliptique $E/K$, la hauteur $h$ v{\'e}rifie
  : \\
\noindent (i) \ \ \ \  $\forall P\in E(\overline{K}) \ \ \ \ \ \ \ \ \ h([m]P)=m^2h(P)+O(1).$\\
\noindent (ii) \ \ \ $\forall P,Q\in E(\overline{K})\ \ \ \ \ h(P+Q)+h(P-Q)=2h(P)+2h(Q)+O(1).$\\
\noindent (iii)\ \ \ $\forall h>0\ \ \ \ \ \ \  \ \ \ \ \ \ \ \textnormal{l'ensemble } \{P\in
E(\overline{K})/\ h(P)\leq h\}
\textnormal{ est fini.}$\\
\noindent Dans les affirmations pr\'ec\'edentes, la constante $O(1)$ d{\'e}pend de $E$ et $m$, mais pas des
points $P$ et $Q$.
\end{prop}

\vspace{.3cm}

\noindent \`A partir de cette hauteur, on peut en construire une plus jolie : la \textit{hauteur de N\'eron-Tate}, not\'ee $\widehat{h}$. La d\'efinition est la suivante : 

\[\widehat{h}(P)=\lim_{n\to+\infty}\frac{h([2^n]P)}{4^n}.\]

\noindent Les propri\'et\'es classiques de cette hauteur sont r\'esum\'ees dans le th\'eor\`eme suivant.

\vspace{.3cm}

\begin{theo} La hauteur canonique est une forme quadratique positive semi-d\'efinie sur $E(\overline{K})$, telle que 
\[ \textit{ d'une part }\forall P\in E(\overline{K})\ \widehat{h}(P)=h(P)+O(1), \textit{ et d'autre part, }\ \widehat{h}(P)=0 \iff P\in E_{\textnormal{tors}}.\]
\end{theo}

\subsection{Multiplication complexe}

\noindent Soient $K$ un corps de nombres et $E/K$ une courbe elliptique \`a multiplication complexe par l'ordre d'un corps quadratique imaginaire $k$. On note $\mathcal{O}_K$ l'anneau d'entiers de $K$ et pour toute place finie $v$ de $K$ on note $k_v$ le corps r\'esiduel associ\'e \`a $v$. Quitte \`a faire une extension de corps ne d\'ependant que de $E/K$ et quitte \`a prendre une courbe elliptique isog\`ene \`a la courbe de d\'epart, on peut supposer que $K$ contient $k$ et que l'anneau des endomorphismes de $E/K$ est exactement $\mathcal{O}_k$, l'anneau des entiers de $k$. De plus, la courbe est \`a multiplication complexe, donc elle a bonne r\'eduction potentielle. Ainsi, quitte \`a remplacer $K$ par une extension de degr\'e born\'e (en fonction de $E/K$), on peut \'egalement supposer que $E/K$ a bonne r\'eduction en toute place de $K$. On fait toutes ces hypoth\`eses dans la suite.

\vspace{.3cm}

\noindent On fixe un point $P\in E(\overline{K})\backslash E_{\textnormal{tors}}$ et on note $D=[K^{\textnormal{ab}}(P):K^{\textnormal{ab}}]$. On choisit alors une extension $F/K$ ab\'elienne finie, telle que $D=[F(P):F]$, ceci \'etant possible car $K^{\textnormal{ab}}$ est le compositum des extensions ab\'eliennes sur $K$.

\vspace{.3cm}

\noindent  Dans la suite, on fixe un mod\`ele de Weierstrass de $E$ de la forme
\[ Y^2=X^3+a_4X+a_6,\]
\noindent o\`u $a_4$ et $a_6$ sont des \'el\'ements de $K$. Si $\wp$ est la fonction de Weierstrass associ\'ee, la courbe complexe $E(\mathbb{C})$ est param\'etr\'ee par $X=\wp(z)$ et $Y=\wp '(z)$.  On rappelle que les points complexes d'une courbe elliptique sont param\'etr\'es par l'isomorphisme de groupes de Lie complexes
\[\mathbb{C}/\Lambda\rightarrow E(\mathbb{C})\ :\  y^2=x^3+a_4x+a_6, \ \ \ z\mapsto \left(\wp(z),\wp'(z)\right),\]
\noindent o\`u $\wp(z)$ est la fonction de Weierstrass d\'efinie par la formule
\[\forall z\in \mathbb{C},\ \ \wp(z)=\frac{1}{z^2}+\sum_{\omega\in \Lambda^*}\left(\frac{1}{(z-\omega)^2}-\frac{1}{\omega^2}\right).\]

\vspace{.3cm}

\noindent Soient $p$ un nombre premier et $v$ une place de $K$ au-dessus de $p$. On rappelle le th{\'e}or{\`e}me fondamental, d\^u \`a Deuring \cite{deuring}, concernant la multiplication complexe que l'on va utiliser ici. On renvoie par exemple {\`a} \cite{si2} Chapter II pour les d\'emonstrations. 

\vspace{.3cm}

\begin{prop} Soit $E/\mathbb{C}$ une courbe elliptique {\`a} multiplication complexe par $\mathcal{O}_k$, anneau d'entiers d'un corps de nombres quadratique imaginaire. Il existe un unique isomorphisme 
$$[\cdot]\ :\ \mathcal{O}_k\rightarrow \textnormal{End}(E)$$
\noindent tel que pour toute diff{\'e}rentielle invariante de $E$, $\omega\in \Omega_E$ et pour tout $\alpha\in \mathcal{O}_k$, on a 
$$[\alpha]^*\omega=\alpha\omega.$$
\noindent De plus le degr{\'e} de $[\alpha]$ est {\'e}gal {\`a} $N_{\mathbb{Q}}^k(\alpha).$
\end{prop}
\vspace{.3cm}

\begin{theo}\textnormal{\textbf{(Deuring)}} Soit $E/K$ une courbe d{\'e}finie sur le corps de nombres $K$, \`a multiplication complexe par le corps quadratique imaginaire $k.$ Soient $p$ un nombre premier et $v$ une place de $K$ au-dessus de $p$ telle que $E$ a bonne r{\'e}duction en $v$. Alors il existe un unique $\alpha_v\in \mathcal{O}_k$ tel que $\widetilde{[\alpha_v]}=\textnormal{Frob}_{E_v}$ o{\`u} $E_v$ est la r{\'e}duction de $E$ sur $\mathbb{F}_v$.
\end{theo}
\vspace{.3cm}

\noindent De plus au cours de la preuve de ce th\'eor\`eme on montre que $q:=N_{\mathbb{Q}}^K(v)=N_{\mathbb{Q}}^K(\alpha)$. Dans les deux lemmes qui suivent, cons{\'e}quences du th{\'e}or{\`e}me pr{\'e}c{\'e}dent, on note $\pi$ une uniformisante dans $k$ de l'id{\'e}al maximal $\mathfrak{M}$ correspondant {\`a} la place $v$.

\vspace{.3cm}

\begin{lemme}\label{fermat} Pour tout {\'e}l{\'e}ment $\alpha\in \mathcal{O}_k$, il existe deux polyn{\^o}mes $R_{\alpha}$ et $S_{\alpha}$ premiers entre eux, {\`a} coefficients dans $\mathcal{O}_k$ tels que 
\[\wp(\alpha z)=\frac{R_{\alpha}(\wp(z))}{S_{\alpha}(\wp(z))},\text{ et }\ \widetilde{S_{\alpha}}\not=0.\]
\noindent Ces deux polyn{\^o}mes sont d{\'e}finis {\`a} multiplication par une m{\^e}me unit{\'e} de $\mathcal{O}_k$ pr{\`e}s. Notamment quand $\alpha=\alpha_v$, on a 
\[R_{\alpha}(X)=uX^q+\pi V(X),\text{ et }\ S_{\alpha}(X)=u+\pi W(X),\]
\noindent o\`u $u$ est une unit\'e $v$-adique de $\mathcal{O}_k$ et $V$ et $W$ sont deux polyn\^omes \`a coefficients dans $\mathcal{O}_k$.
\end{lemme}

\vspace{.3cm}

\begin{lemme} \label{premier} Pour tout $\alpha\in \mathcal{O}_k$, les polyn{\^o}mes $\widetilde{R_{\alpha}}$ et $\widetilde{S_{\alpha}}$ sont premiers entre eux.
\end{lemme}
\demo On trouvera par exemple une preuve de ces deux lemmes dans \cite{laurent} lemmes 3.1 et 3.2 respectivement.\hfill$\Box$

\vspace{.3cm}

\noindent Nous aurons \'egalement besoin d'un lemme sur les endomorphismes du groupe formel associ\'e \`a la courbe elliptique $E$. Si $P$ est un point de la courbe de coordonn\'ees affines ($X,Y)$, on note $t=-\frac{X}{Y}$ et on note $[\alpha_v]$ l'op\'erateur du groupe formel associ\'e \`a l'endomorphisme $\alpha_v$.

\vspace{.3cm}

\begin{lemme}\label{prep} Il existe une s\'erie enti\`ere $\psi$, \`a coefficients dans $\mathcal{O}_{K_v}$ telle que 
\[ [\alpha_v](t)=t^p+\pi_p\psi(t).\]
\end{lemme}
\demo C'est le lemme 3.3 de \cite{laurent}.\hfill$\Box$

\section{R\'eductions}\label{reduc}

\noindent On fait maintenant les m\^emes r\'eductions que dans le cas multiplicatif d\^u \`a Amoroso-Zannier. On note $\mathcal{P}$ l'ensemble des nombres premiers qui se d\'ecomposent totalement dans $K$. Pour chacune des places $v$ de $K$ au-dessus d'un tel premier $p$, la compl\'etion $v$-adique de $K$ est $K_v=\mathbb{Q}_p$. Pour $p\in \mathcal{P}$, on notera donc $K_p$ cette compl\'etion dans la suite. Soient $p\in \mathcal{P}$ et $F/K$ une extension ab\'elienne finie, on note $e_p(F)$ l'indice de ramification de $p$ dans $F$ et $F_v$ la compl\'etion $v$-adique de $F$ en $v$. On a $K_p=\mathbb{Q}_p$, donc $F_v$ est une extension ab\'elienne de $\mathbb{Q}_p$. Par le th\'eor\`eme de Kronecker-Weber local, elle est donc contenue dans une extension cyclotomique de $\mathbb{Q}_p$ que l'on notera $\mathbb{Q}_p(\zeta_m)$. On pose $m=m_p(F)$ le plus petit entier ayant cette propri\'et\'e et on d\'efinit $\mathfrak{f}_p(F)$ le \textit{conducteur local de $F$ en $p$} comme \'etant la plus grande puissance de $p$ divisant $m$ (il s'agit bien du conducteur local au sens de la th\'eorie du corps de classes local). On pose 
\[\mathfrak{f}(F)=\prod_{p\in\mathcal{P}}\mathfrak{f}_p(F),\]
\noindent le \textit{conducteur de $F$} et on note que si $F'\subset F$ alors $\mathfrak{f}(F')\leq\mathfrak{f}(F)$.

\vspace{.3cm}

\noindent Soit maintenant $P$ un point de $E(\overline{K})\backslash E_{\textnormal{tors}}$ contredisant le th\'eor\`eme \ref{thmoi}, de degr\'e minimal, \textit{i.e.}, tel que pour tout point $P'\in E(\overline{K})\backslash E_{\textnormal{tors}}$ de degr\'e $D'<D$ sur $K^{\textnormal{ab}}$, on a 
\[\widehat{h}(P')\geq\frac{c(E/K)}{D'}\left(\frac{\log \log 5D'}{\log 2D'}\right)^{13}.\]

\vspace{.3cm}

\begin{lemme}\label{r1} Pour d\'emontrer le th\'eor\`eme \ref{thmoi}, on peut supposer que pour tout point de torsion $T\in E_{\textnormal{tors}}$ on a  $[K^{\textnormal{ab}}(P+T):K^{\textnormal{ab}}]\geq D.$
\end{lemme}
\demo La hauteur de N\'eron-Tate est invariante par translation par un point de torsion. Le r\'esultat d\'ecoule donc imm\'ediatement de la d\'efinition du point $P$ et de la d\'ecroissance pour $t\geq 1$ de la fonction $t\mapsto \frac{c(E/K)}{t}\left(\frac{\log \log 5t}{\log 2t}\right)^{13}$.\hfill$\Box$

\vspace{.3cm}

\noindent Soit $\mathcal{A}$ l'ensemble des extensions ab\'eliennes finies $F/K$ telles qu'il existe un point de torsion $T\in E_{\textnormal{tors}}$ tel que $[F(P+T):F]\leq D$, \textit{i.e.}, tel que $[F(P+T):F]=D$ par le lemme pr\'ec\'edent. Cet ensemble est non vide, puisque par d\'efinition de $K^{\textnormal{ab}}$, on sait qu'il existe une extension ab\'elienne finie $F/K$ telle que $[F(P):F]=[K^{\textnormal{ab}}(P):K^{\textnormal{ab}}]=D.$ L'extension $F$ et le point $T=0$ montrent donc que $\mathcal{A}$ est non vide. On d\'efinit alors l'entier
\[\mathfrak{f}=\min_{F\in\mathcal{A}}\mathfrak{f}(F).\]

\vspace{.3cm}

\begin{lemme}\label{reduction}Avec les notations pr\'ec\'edentes, pour d\'emontrer le th\'eor\`eme \ref{thmoi}, on peut supposer que 
\begin{equation}\label{hyp1}
D=[F(P):F] \textnormal{ o\`u $F/K$ est une extension appartenant \`a $\mathcal{A}$, contenue dans $K(P)$.}
\end{equation}
\noindent On peut \'egalement supposer que 
\begin{equation}\label{hyp2}
\mathfrak{f}(F)=\mathfrak{f}.
\end{equation}
\noindent Enfin, on peut aussi supposer que  
\begin{equation}\label{hypocorps} 
\forall T\in E_{\textnormal{tors}} \textnormal{ tel que }K(P+T)\subset K(P), \textnormal{ on a } K(P+T)=K(P).
\end{equation}
\end{lemme}
\demo Par d\'efinition de $K^{\textnormal{ab}}$, il existe une extension ab\'elienne finie $F/K$ telle que $D=[F(P):F]$, donc appartenant \`a $\mathcal{A}$. On prend dans $\mathcal{A}$ une extension $F/K$ r\'ealisant le min des $\mathfrak{f}(F)$, \textit{i.e.}, r\'ealisant $\mathfrak{f}$. Montrons qu'on peut supposer (\ref{hypocorps}). Soit $T\in E_{\textnormal{tors}}$ tel que $K(P+T)\subset K(P)$, alors, le point $T$ est d\'efini sur le corps de nombres $K(P)$ car $P+T-P=T$. Si pour tous ces $T$, on a l'\'egalit\'e $K(P+T)=K(P)$, il n'y a alors rien \`a montrer. Sinon, on note $\mathcal{T}$ l'ensemble fini des points de torsion tels que $K(P+T)\subsetneq K(P)$. Soient $T\in \mathcal{T}$ et $P_1=P+T$. L'extension $K(P_1)$ est une sous-extension stricte de $K(P)$. On note  $\mathcal{T}_1$ l'ensemble fini des points de torsion tels que $K(P_1+T)\subsetneq K(P_1)$. Si $\mathcal{T}_1$ est non vide, on choisit $T_1\in\mathcal{T}_1$ et on pose $P_2=P_1+T_1$. L'extension $K(P_2)$ est une sous-extension stricte de $K(P_1)$. On construit ainsi une cha\^ine 
\[K(P_n)\subsetneq\ldots\subsetneq K(P_1)\subsetneq K(P).\]
\noindent Donc pour $n$ assez grand, on sait que $K(P_{n+1})=K(P_n)$, autrement dit que l'ensemble $\mathcal{T}_{n+1}$ correspondant est vide, c'est-\`a-dire que 
\[\forall T\in E_{\textnormal{tors}} \textnormal{ tel que }K(P_n+T)\subset K(P_n), \textnormal{ on a } K(P_n+T)=K(P_n).\]
\noindent Or par construction, on a $P_n=P+T_n$ o\`u $T_n$ est un point de torsion de $E$, donc le lemme \ref{r1} pr\'ec\'edent assure que $[F(P_n):F]\geq D_n\geq D$. De plus on a 
\[D_n=[K^{\textnormal{ab}}(P_n):K^{\textnormal{ab}}]\leq [F(P_n):F]\leq [F(P):F]=D\]
\noindent car $K(P_n)\subset K(P)$, donc $D_n=D$. La hauteur de N\'eron-Tate \'etant invariante par translation par un point de torsion, on a \'egalement $\widehat{h}(P_n)=\widehat{h}(P)$. Enfin, quitte \`a remplacer $F$ par $F_1=F\cap K(P_n)$, on voit que l'on peut aussi supposer l'hypoth\`ese (\ref{hyp1}) vraie. La fonction $\mathfrak{f}(\cdot)$ \'etant croissante, l'hypoth\`ese (\ref{hyp2}) est elle aussi v\'erifi\'ee, ce qui conclut.\hfill$\Box$

\vspace{.3cm}

\noindent Dans toute la suite on supposera d\'esormais vraies les hypoth\`eses (\ref{hyp1}), (\ref{hyp2}) et (\ref{hypocorps}). 

\vspace{.3cm}

\begin{rem} \label{rem1}\textnormal{On note que, comme $K\subset F\subset K(P)$, on a aussi, $F(P)=K(P)$.}
\end{rem}

\vspace{.3cm}

\noindent On peut maintenant \'enoncer les deux lemmes de r\'eduction qui nous serviront dans la suite. Le premier est inspir\'e du Lemma 2.1. (ii) de \cite{AZ}, le second est plus classique dans le cadre du prob\`eme de Lehmer.

\vspace{.3cm}

\begin{lemme}\label{K}Soient $p\in \mathcal{P}$ et $v$ une place de $K$ au-dessus de $p$, alors, pour d\'emontrer le th\'eor\`eme \ref{thmoi}, on peut supposer que 
\[ \textit{soit }K(\alpha_v(P))=K(P),\textit{ soit }[K(P):K(\alpha_v(P))]=p.\]
\end{lemme}
\demo On consid\`ere le diagramme
\[
\xymatrix{
 		& K\left(P,E[\alpha_v]\right) 	&  							\\
K(P)	\ar[ur]	&   				& K\left(\alpha_v(P),E[\alpha_v]\right)\ar[ul]_{\textnormal{\large{G}}}	\\
		& K(\alpha_v(P))\ar[ul]\ar[ur]	&							
}\]
\noindent L'extension $K(P,E[\alpha_v])/K(\alpha_v(P),E[\alpha_v])$ est galoisienne d'ordre $1$ ou $p$. En effet on a une injection naturelle $\textnormal{Gal}(\overline{K}/K(\alpha_v(P),E[\alpha_v]))\hookrightarrow E[\alpha_v]$ : les conjugu\'es de $P$ par l'action de $\textnormal{Gal}(\overline{K}/K(\alpha_v(P),E[\alpha_v]))$ sont parmi les $P+T$, o\`u $T\in E[\alpha_v]$ et $\alpha_v$ est une isog\'enie cyclique d'ordre $p$. 

	Si le groupe de Galois correspondant G est d'ordre $p$, alors l'extension $K(P)/K(\alpha_v(P))$ est \'egalement d'ordre $p$. 

	Si G est d'ordre $1$, on va montrer qu'il existe $T\in E[\alpha_v]$ tel que $K(P+T)\subset K(\alpha_v(P))$. On regarde l'action de $\textnormal{Gal}(\overline{K}/K(\alpha_v(P))$ sur l'ensemble $\left\{P+T\ /\ T\in \textnormal{ker}[\alpha_v]\right\}$. Soit il y a une seule orbite, auquel cas $[K(\alpha_v(P)):K(P)]=p$ ; soit l'orbite $\omega_P$, contenant $P$ est de cardinal $m$ strictement inf\'erieur \`a $p$, donc premier \`a $p$. Dans ce cas, il existe $T'\in E[\alpha_v]$, tel que
\[\sum_{T\in \omega_P}(P+T)=mP+T'\]
\noindent est stable sous l'action de $\textnormal{Gal}(\overline{K}/K(\alpha_v(P))$. Par le th\'eor\`eme de B\'ezout, il existe deux entiers, $\lambda$ et $\mu$ tels que $\lambda m+\mu p=1$. Par ailleurs, en notant $\alpha_v^{\vee}$ l'isog\'enie duale de $\alpha_v$, on a 
\[ K\left([p]P\right)=K\left(\alpha_v^{\vee}(\alpha_v(P))\right)\subset K(\alpha_v(P)).\]
\noindent Ainsi, on a les inclusions
\[K\left([\lambda]([m]P+T')+[\mu p]P)\right)=K\left(P+[\lambda]T'\right)\subset K(\alpha_v (P)).\]
\noindent On a donc l'inclusion $K(P+[\lambda]T')\subset K(P)$. Par l'hypoth\`ese (\ref{hypocorps}) ceci entra\^ine que 
\[K(P)=K(P+[\lambda]T')\subset K(\alpha_v(P)).\]
\noindent On en d\'eduit que $K(P)=K(\alpha_v(P)).$ \hfill$\Box$

\vspace{.3cm}

\begin{lemme}\label{L} Pour tout $p\in \mathcal{P}$ sauf pour au plus $\frac{1}{2}\log D$ d'entre eux et pour toute place $v$ de $K$ au-dessus de $p$, on a 
\[F(P)=F(\alpha_v(P)).\]
\end{lemme}
\demo C'est le lemme combinatoire classique de Dobrowolski \cite{dob} (d\^u \`a Laurent \cite{laurent} lemme 4.2 dans le cas des courbes elliptiques).\hfill$\Box$

\vspace{.3cm}

\noindent Notons que l'on pourrait \'eviter de recourir \`a ce lemme combinatoire, en faisant un raisonnement du m\^eme type que dans le lemme pr\'ec\'edent, comme il est fait dans l'article d'Amoroso-Zannier \cite{AZ}. Dans la suite on notera $\mathcal{P}^{\star}$ le sous-ensemble de $\mathcal{P}$ form\'e des premiers  v\'erifiant le lemme \ref{L} pr\'ec\'edent.

\section{Lemmes d'extrapolation}
\noindent Dans la partie \ref{extra} on va faire deux extrapolations diff\'erentes, selon qu'il y a beaucoup de places de $F$ au-dessus de $\mathcal{P}^{\star}$ ayant un gros indice de ramification, ou non, tout ceci \'etant bien \'evidemment quantifi\'e. On commence par les lemmes qui nous permettront d'extrapoler dans le cas o\`u il y a  beaucoup de ramification.

\subsection{Lemme ramifi\'e}

\begin{lemme}\label{baker}Soient $E/K$ une courbe elliptique \`a multiplication complexe, $v$ une place de bonne r\'eduction ordinaire et $I_v$ le groupe d'inertie de $\textnormal{Gal}(\overline{K_v}/K_v)$. Alors, pour tout entier $n\geq 1$, on a l'isomorphisme de $I_v$-module $E[\alpha_v^n]\simeq \mu_{p^n}$.
\end{lemme}
\demo C'est le lemme 3.2 de \cite{baker}.\hfill$\Box$

\vspace{.3cm}

\noindent Le lemme suivant est inspir\'e du lemme 3.2. de \cite{AZ}.

\vspace{.3cm}

\begin{lemme}\label{ramifie} Soient $p\in \mathcal{P}^{\star}$ et $e_p$ son indice de ramification dans $F$. Il existe un sous-groupe $H_p$ de $\textnormal{Gal}(F/K)$ d'ordre 
\[\mid H_p\mid\geq\min\{e_p,p\},\]
\noindent tel que 
\[\mid x^p-\sigma x^p\mid_w\leq \frac{1}{p},\]
\noindent pour tout $x\in \mathcal{O}_F$, tout $\sigma\in H_p$ et toute place $w$ de $F$ au-dessus de $p$. De plus, pour toute place $v$ de $K$ au-dessus de $p$ et pour toute extension $\tau\in \textnormal{Gal}(\overline{K}/K)$ de $\sigma\in H_p-\{\textnormal{Id}\}$, on a 
\[\tau(\alpha_v(P))\not=\alpha_v(P).\]
\end{lemme}
\demo La fabrication de $H_p$ et l'estimation de son cardinal se fait comme dans l'article de Amoroso-Zannier : soient $v$ une place de $F$ \'etendant $p$ et $F_v$ le compl\'et\'e $v$-adique de $F$. On pose $m:=m_p(F)$ le plus petit entier $m$ tel que $F_v\subset\mathbb{Q}_p(\zeta_m).$ On d\'ecompose $m$ sous la forme $m=\mathfrak{f}_p\cdot n$ o\`u $n$ est premier \`a $p$ et $\mathfrak{f}_p$ est le conducteur local de $F$ en $p$.

Si $p$ ne se ramifie pas dans $F$, alors $e_p=1$ et $H_p=\{\textnormal{Id}\}$ convient. On peut donc supposer que $p$ se ramifie dans $F$, donc \textit{a fortiori} dans $\mathbb{Q}(\zeta_m)$. Ainsi $p$ divise le conducteur local $\mathfrak{f}_p$. On pose $\Sigma_p$ le groupe de Galois de l'extension $\mathbb{Q}_p(\zeta_m)/\mathbb{Q}_p(\zeta_{m/p})$. C'est un groupe cyclique d'ordre $p$ ou $p-1$ selon que $p^2$ divise $\mathfrak{f}_p$ ou non. Par la propri\'et\'e de minimalit\'e de $m$, $\Sigma_p$ ne fixe pas $F_v$, donc induit par restriction un sous-groupe non-trivial $H_v^*$ de $\textnormal{Gal}(F_v/K_p)$. On note que si $p^2$ ne divise pas $\mathfrak{f}_p$, alors l'ordre de $H_v^*$ est au moins $e_p$ car l'extension $\mathbb{Q}_p(\zeta_{m/p})/\mathbb{Q}_p$ est non-ramifi\'ee ; alors que si $p^2$ divise $\mathfrak{f}_p$, n\'ecessairement $H_v^*$ est d'ordre $p$. On d\'efinit $H_v$ comme \'etant l'image isomorphe de $H_v^*$ dans $\textnormal{Gal}(F/K)$. On peut voir que $H_v^*$ ne d\'epend pas de $v$, mais seulement de $p$. Il en est de m\^eme de $H_v$ que l'on note d\'esormais $H_p$. On a d\'eja obtenu l'estimation de son cardinal. 

Montrons la propri\'et\'e de congruence : soit $\mathcal{O}$ l'anneau des entiers de $\mathbb{Q}_p(\zeta_m)$. On a 
\begin{equation}\label{cyclo} 
\forall x\in\mathcal{O},\ \forall \sigma\in \Sigma_p \ \ x^p\equiv\sigma x^p \mod p\mathcal{O},
\end{equation}

\noindent (cf. par exemple \cite{AZ} p. 717). Ainsi pour tout $x\in\mathcal{O}_F$ et pour tout $\sigma\in H_p$, l'entier $x^p-\sigma x^p\in F$ est d'ordre sup\'erieur \`a $e_p$ en $v$. 

Montrons maintenant la derni\`ere propri\'et\'e. Soient $\sigma\in H_p-\{\textnormal{Id}\}$ et $\tau\in\textnormal{Gal}(\overline{K}/K)$ une extension de $\sigma$. Supposons par l'absurde que $\tau(\alpha_v(P))=\alpha_v(P)$.  Soit $\mathbb{E}$ le sous-corps de $F$ fixe par $\sigma$. On a $[\mathbb{E}(\alpha_v(P)):\mathbb{E}]=[F(\alpha_v(P)):F]$. De plus, par le lemme \ref{L}, le point $P$ est d\'efini sur la m\^eme extension de $F$ que $\alpha_v(P)$. Donc, 
\begin{equation}\label{D}
[\mathbb{E}(\alpha_v(P)):\mathbb{E}]=[F(\alpha_v(P)):F]=[F(P):F]=D.
\end{equation}

\vspace{.3cm}

\noindent On va maintenant montrer que $\mid \Sigma_p\mid=p$ : tout d'abord, comme $\mathbb{E}$ est strictement inclus dans $F$, on a $[F(P):\mathbb{E}]>[F(P):F]$ et donc, d'apr\`es (\ref{D}), 
\begin{equation}\label{tri}
 [F(P) :E(\alpha_v(P))]=\frac{[F(P):E]}{[E(\alpha_v(P)):\mathbb{E}]}=\frac{[F(P):E]}{[F(P):F]}>1.
\end{equation}
\noindent Par ailleurs, d'apr\`es la remarque \ref{rem1}, on a $K(P)=F(P)$ et $K\subset \mathbb{E}(\alpha_v(P))\subset F(P)$. Ainsi, $[F(P):\mathbb{E}(\alpha_v(P))]$ divise $[K(P):K(\alpha_v(P))]$. Si $K(P)=K(\alpha_v(P))$, alors $\tau$ fixe $K(P)=F(P)$ donc fixe $F$ ce qui contredit le choix de $\sigma\not=\textnormal{Id}$. Ainsi, par le lemme \ref{K}, l'extension $K(P)/K(\alpha_v(P))$ est de degr\'e $p$. On en d\'eduit que l'extension $F(P)/\mathbb{E}(\alpha_v(P))$ qui est non triviale par (\ref{tri}), est de degr\'e $p$. On a ainsi 
\[ [F(\alpha_v(P):\mathbb{E}(\alpha_v(P))]=\frac{[F(P):\mathbb{E}(\alpha_v(P))]}{[F(P):F(\alpha_v(P))]}=[F(P):\mathbb{E}(\alpha_v(P))]=p \textnormal{ par le lemme \ref{L}}.\]
\noindent L'extension $F/\mathbb{E}$ \'etant galoisienne, on en d\'eduit que $\mid H_p\mid=[F:E]=p$, donc par construction de $H_p$ on obtient $\mid \Sigma_p\mid=p$.

\vspace{.3cm}

	On sait que sur une courbe elliptique \`a multiplication complexe, on a bonne r\'eduction ordinaire en toutes les places $v$ au-dessus d'un premier $p\in\mathcal{P}.$ On peut donc appliquer le lemme \ref{baker} dans notre situation. Par ce lemme on sait que les points de $\alpha_v^k$-torsion sont d\'efinis sur $\mathbb{Q}_p(\zeta_{p^k})\subset\mathbb{Q}_p(\zeta_{m})$. Ainsi $F_v(E[\alpha_v^k])\subset \mathbb{Q}(\zeta_m)$, donc le groupe de Galois $\Sigma_p$ induit par restriction un sous-groupe non-trivial de $\textnormal{Gal}(F(E[\alpha_v^k])/K)$ qui est cyclique d'ordre $p$ par le paragraphe pr\'ec\'edent. Soit donc $\mathbb{F}\subset F(E[\alpha_v^k])$ son sous-corps fixe. Soient $x\in \mathbb{E}$ et $\rho \in G_1=\textnormal{Gal}(F(E[\alpha_v^k])/\mathbb{F})-\{\textnormal{Id}\}$. Puisque $[F:\mathbb{E}]=p$, le morphisme $\sigma$ engendre le groupe $\textnormal{Gal}(F/\mathbb{E}),$ donc il existe un entier $u$ tel que $\rho_F=\sigma^u$. Notamment, on en d\'eduit que $\rho(x)=x$, c'est-\`a-dire que 
\begin{equation}\label{EinF}
 \mathbb{E}\subset \mathbb{F}.
\end{equation}

\vspace{.3cm}

	On va maintenant montrer qu'il existe un point de $\alpha_v^k$-torsion $T$ tel qu'on ait l'inclusion $\mathbb{F}(P+T)\subset \mathbb{F}(\alpha_v(P))$. Si $\mathbb{F}(P)\subset\mathbb{F}(\alpha_v(P))$, il n'y a rien \`a montrer. Sinon, on a \textit{a fortiori} l'inclusion stricte 
\[\mathbb{F}(\alpha_v(P))\subsetneq F\left(E[\alpha_v^k],P\right).\]
\noindent Les extensions \'etant galoisiennes, $[F(E[\alpha_v^k],\alpha_v(P)):\mathbb{F}(\alpha_v(P))]$ divise $[F(E[\alpha_v^k]):\mathbb{F}]=p$. De plus, par le lemme \ref{L}, $F(E[\alpha_v^k],P)=F(E[\alpha_v^k],\alpha_v(P))$, donc on a l'\'egalit\'e 
\[ \left[F\left(E[\alpha_v^k],P\right):\mathbb{F}\left(\alpha_v(P)\right)\right]=p.\]
\noindent Ainsi, le morphisme de restriction 
\[\textnormal{ res : } \textnormal{Gal}\left(F(E[\alpha_v^k],P)/\mathbb{F}(\alpha_v(P))\right)\rightarrow \textnormal{Gal}\left(F(E[\alpha_v^k])/\mathbb{F}\right),\]
\noindent entre groupes de m\^eme cardinaux est un isomorphisme. Soit $\widetilde{\rho}$ un g\'en\'erateur du groupe cyclique $\textnormal{Gal}\left(F(E[\alpha_v^k],P)/\mathbb{F}(\alpha_v(P))\right)$... Il existe un point de $\alpha_v$-torsion $T_1$ tel que 
\[ \widetilde{\rho}(P)=P+T_1.\]
\noindent De plus, par le lemme \ref{baker}, on a l'isomorphisme de $I_v$-modules, $E[\alpha_v^k]\simeq \mu_{p^k}$, donc si $T_2$ est un point de $E[\alpha_v^k]\backslash E[\alpha_v^{k-1}]$, alors le point $T_3=\rho(T_2)-T_2$ est d'ordre $p$. Finalement, il existe un entier $v$ tel que 
\[T_1=vT_3.\]
\noindent On pose $T=-rT_2$. On a alors
\[\rho(P+T)=P+T_1-v\rho(T_2)=P+vT_3-vT_3-vT_2=P+T.\]
\noindent Ceci nous donne bien l'inclusion $\mathbb{F}(P+T)\subset \mathbb{F}(\alpha_v(P))$. 

\vspace{.3cm}

	En utilisant (\ref{D}) et (\ref{EinF}), on obtient
\[\left[\mathbb{F}(P+T):\mathbb{F}\right]\leq \left[\mathbb{F}(\alpha_v(P):\mathbb{F}\right]\leq\left[\mathbb{E}(\alpha_v(P):\mathbb{E}\right]=D.\]
\noindent Or par construction, $\mathbb{F}_v\subset \mathbb{Q}_p(\zeta_{m/p})$ et $\mathbb{F}\subset F(E[\alpha_v^k])$, donc
\[ \mathfrak{f}_p(\mathbb{F})\leq \frac{p^k}{p}< \mathfrak{f}_p(F),\textnormal{ et, si $l\not=p,$ }  \mathfrak{f}_l(\mathbb{F})\leq \mathfrak{f}_l(F).\]
\noindent On en conclut, que 
\[\left[\mathbb{F}(P+T):\mathbb{F}\right]\leq D\textnormal{ et, } \mathfrak{f}(\mathbb{F})< \mathfrak{f}(F)=\mathfrak{f},\]
\noindent ce qui contredit la d\'efinition de $\mathfrak{f}$. Ceci conclut la preuve par l'absurde.\hfill $\Box$

\subsection{Lemme non-ramifi\'e}
\noindent On passe maintenant au lemme qui va nous permettre de faire l'extrapolation dans le cas o\`u il n'y a pas beaucoup de ramification. Il s'agit du m\^eme lemme que dans le cas multiplicatif.

\vspace{.3cm}

\begin{lemme}\label{nonramifie} Soit $p\in\mathcal{P}^{\star}$, il existe $\Phi_p\in\textnormal{Gal}(F/K)$ tel que
\[\mid x^p-\Phi_p x\mid_v\leq p^{-\frac{1}{e_p}},\]
\noindent o\`u $x\in \mathcal{O}_F$ et $v$ est une valuation sur $\overline{\mathbb{Q}}$ \'etendant $p$.
\end{lemme}
\demo C'est le lemme 3.1. de \cite{AZ}.\hfill$\Box$

\vspace{.3cm}

\section{Lemme de Siegel}\label{sieg}
\noindent Dans la suite, on consid\`ere un point $P_1$ qui sera soit $P$ soit $\alpha_v(P)$. On note $\wp(u)$ la coordonn\'ee $x$ de $P_1$ et on note $\wp(u_1),\ldots,\wp(u_D)$ les diff\'erents conjugu\'es de $\wp(u)$ sur $F$. On dit que les $u_i$ sont \textit{les conjugu\'es de $u$}.

\vspace{.3cm}

\noindent Soient $L$ et $T$ deux entiers strictement positifs et $N\in ]\sqrt{L},2\sqrt{L}[$ un nombre premier (qui existe par le ``postulat de Bertrand''). On va construire une fonction
\[ \varphi(z)=\sum_{\lambda_1=0}^{L}\sum_{\lambda_2=0}^{L}p(\lambda_1,\lambda_2)\wp(z)^{\lambda_1}\wp(Nz)^{\lambda_2}\]
\noindent avec $p(\lambda_1,\lambda_2)\in \overline{\mathbb{Z}}$ non tous nuls, telle que  : $\varphi$ n'est pas constamment nulle sur $E(\mathbb{C})$, $\varphi$ est nulle en les conjugu\'es $u_i$ de $u$ avec multiplicit\'e $T$ et les coefficients $p(\cdot,\cdot)$ sont bien control\'es. Le premier point est assur\'e par le choix de $N$ et le fait que les $p(\cdot,\cdot)$ ne sont pas tous nuls, le second d\'ecoule d'un lemme de Siegel absolu.

\vspace{.3cm}

\begin{lemme}\label{siegelabs} Soient $n$ un entier et $S$ un sous $\overline{\mathbb{Q}}$-espace vectoriel de dimension $d$ de $\overline{\mathbb{Q}}^n$. Pour tout $\varepsilon>0$, il existe un vecteur $\textbf{x}\in S$ tel que 
\[h_2(\textbf{x})\leq \frac{h_2(S)}{d}+\frac{\log d}{2}+\varepsilon.\]
\end{lemme}
\demo cf. \cite{davphi2} lemme 4.7 et la remarque qui suit.\hfill$\Box$

\vspace{.3cm}

\begin{prop}\label{siegelth} Soient $L$, $T$ et $k$ trois entiers positifs tels que $L^2\geq kT$ et $k^{c_2}>(L+T)^2$ pour une certaine constante absolue $c_2>0$. Avec les notations pr\'ec\'edentes, on peut construire la fonction $\varphi$, s'annulant en $u_1,\ldots,u_k$ avec multiplicit\'e sup\'erieure \`a $T$, telle que
\[h_2(\varphi)\leq  \frac{ckT}{(L+1)^2-kT}\left(LN^2\widehat{h}(P_1)+T\log(T+L)+T\log N+L\right)+\log L,\]
\noindent o\`u $c$ est une constante ne d\'ependant que de $E/K$.
\end{prop}
\demo  Par r{\'e}currence sur $t\leq T$, on montre qu'il existe un polyn{\^o}me $Q_{\lambda_1,\lambda_2,t}$ dans $\mathcal{O}_k[X_1,\ldots,X_4]$ de degr{\'e} partiel en chaque variable major{\'e} par $L+2t$, {\`a} coefficients de valeur absolue major{\'e}e par $c_1k^{c_2t}$ et tel que 
$$\frac{d^t}{dz^t}\left(\wp(z)^{\lambda_1}\wp(Nz)^{\lambda_2}\right)=Q_{\lambda_1,\lambda_2,t}\left(\wp(z),\wp'(z),\wp(Nz),\wp'(Nz)\right).$$

\noindent Au rang $t=0$, le polyn{\^o}me $Q=X_1^{\lambda_1}X_3^{\lambda_2}$ convient.

\noindent Supposons la propri\'et\'e vraie au rang $t$ et montrons-la au rang $t+1$ : en notant abusivement $Q_t$ le polyn\^ome $Q_{\lambda_1,\lambda_2,t}$, on a 
\begin{align*}
\frac{d^{t+1}}{dz^{t+1}}\left(\wp(z)^{\lambda_1}\wp(Nz)^{\lambda_2}\right)    & = \frac{d}{dz}\frac{d^t}{dz^t}\left(\wp(z)^{\lambda_1}\wp(Nz)^{\lambda_2}\right)\\
                                & = \frac{d}{dz} Q_t\left(\wp(z),\wp'(z),\wp(Nz),\wp'(Nz)\right)\\
                                & = \frac{\partial Q_t}{\partial X_1}\left(\cdot\right)\wp'(z)+\ldots+N\frac{\partial Q}{\partial X_4}\left(\cdot\right)\wp''(Nz).\\
\end{align*}
\noindent En utilisant la relation $\wp''(Nz)=6\wp(Nz)^2+2a_4$, on pose donc 
\[Q_{t+1}= X_2\frac{\partial Q_t}{\partial X_1}+\ldots+N(6X_3^2+2a_4)\frac{\partial Q_t}{\partial X_4}.\]
\noindent On a clairement $\textnormal{deg}_{X_i}Q_{t+1}\leq L+2(t+1)$. De plus, en notant $q_{i,t}$ les coefficients de $Q_t$, on a 
\[\vert q_{i,t+1}\vert \leq 6c_1N(L+2t)k^{c_2t}\leq 12c_1\sqrt{L}(L+2T)k^{c_2t}\leq c_1k^{c_2(t+1)}.\]
\noindent Finalement, le syst{\`e}me $\forall t\leq T-1$, $\forall i\in[\![1,k]\!]$, $\varphi^{(t)}(u_i)=0$ s'{\'e}crit : 
\[ \forall  t\leq T-1,\ \  \forall i\in[\![1,k]\!],\ \  \sum_{\lambda_1=0}^{L}\sum_{\lambda_2=0}^{L} p(\lambda_1,\lambda_2)Q_{\lambda_1,\lambda_2,t}\left(\wp(u_i),\wp'(u_i),\wp(Nu_i),\wp'(Nu_i)\right)=0.\]

\noindent Pour tout $0\leq i\leq k$, posons 
\[\forall 0\leq t\leq T-1,\  \forall 0\leq \lambda_1,\lambda_2\leq L,\ \ \alpha_{(\lambda_1,\lambda_2),t}^{(i)}=Q_{\lambda_1,\lambda_2,t}\left(\wp(u_i),\wp'(u_i),\wp(Nu_i),\wp'(Nu_i)\right).\]

\noindent On consid\`ere les vecteurs 
\[\textbf{y}_{i,t}=\left(\alpha_{(0,0),t}^{(i)},\ldots,\alpha_{(L,0),t}^{(i)},\ldots,\alpha_{(L,L),t}^{(i)}\right)\in\overline{\mathbb{Q}}^{(L+1)^2}. \]
\noindent Comme dans \cite{davidhindry} p.42 in\'egalit\'e (18) et suivante, on v\'erifie que les coefficients du syst\`eme
\[  \textbf{y}_{i,t}\cdot \textbf{x}=0,\ \ 0\leq i\leq k,\ \ 0\leq t\leq T-1\]
\noindent avec $\textbf{x}\in \overline{\mathbb{Q}}^{(L+1)^2} $ sont tous de hauteur au plus
\[c_3\left(LN^2\widehat{h}(P_1)+T\log(T+L)+T\log N+L\right).\]
\noindent Par ailleurs, le $\overline{\mathbb{Q}}$-espace vectoriel 
\[S=\left\{\textbf{x}\in\overline{\mathbb{Q}}^{(L+1)^2}\ \ / \ \ \textbf{y}_{i,t}\cdot \textbf{x}=0,\ \ 0\leq i\leq k,\ \ 0\leq t\leq T-1\ \right\}\]
est de dimension $(L+1)^2-kT$ et les vecteurs $\textbf{y}_{i,t}$ forment une base de l'orthogonal $S^{\bot}$.  De plus, le Lemma IV p.10 de \cite{schmidt}, nous indique que 
\[h_2\left(S\right)=h_2\left(S^{\bot}\right)\leq \sum_{i,t}h_2(\textbf{y}_{i,t})\leq c_3kT\left(LN^2\widehat{h}(P_1)+T\log(T+L)+T\log N+L\right).\]
\noindent Pour conclure, on applique le lemme \ref{siegelabs} avec $\varepsilon=\frac{1}{2}\log\frac{(L+1)^2}{(L+1)^2-kT}$.\hfill$\Box$

\section{Extrapolation}\label{extra}
\noindent Il y a deux cas, selon que l'on a ``beaucoup'' de premiers ayant ``beaucoup'' de ramification ou non (ceci \'etant quantifi\'e). On commence par le cas qui sera utilis\'e quand il n'y a pas beaucoup de grande ramification.

\vspace{.3cm}

\begin{prop}\label{nr}Soient $L_1$ et $T_1$ deux entiers strictement positifs d'ordre de grandeur polynomial en $D$, tels que $L_1^2\geq DT_1$. On pose $P_1=P$ et  on consid\`ere la fonction $\varphi$ obtenue dans la proposition \ref{siegelth} avec $L=L_1$, $T=T_1$ et $k=D$. Soient $p\in\mathcal{P}^{\star}$ et $v$ une place \'etendant $p$ sur $\overline{K}$. Pour tout $t\leq L_1$ et pour tout $\tau\in\textnormal{Gal}(\overline{K}/K)$ \'etendant le morphisme $\Phi_p$ du lemme \ref{nonramifie}, on a
\[\log \mid\tau(\varphi)^{(t)})(\alpha_p u)\mid_v\leq -\frac{T_1}{2e_p}\log p+8L_1\log\max\{1,\mid\wp(N\alpha_vu)\mid_v\}.\]
\end{prop}
\demo Il s'agit essentiellement du deuxi\`eme pas de \cite{laurent} \`a la diff\'erence que l'on utilise un lemme de Siegel absolu, ce qui conduit \`a supposer l'annulation en un point et en tous ses conjugu\'es, ainsi qu\`a faire intervenir l'indice de ramification $e_p$ de $p$ dans $L$. On \'etend $v$ au corps $\overline{K}(X)$ en posant $\mid X\mid_v=1$.

\vspace{.3cm}

\noindent Il y a deux cas : soit $\wp(u)$ est un $v$-entier, soit non.

\vspace{.3cm}

\noindent Cas 1 : on v\'erifie simplement que l'on peut \'ecrire $\varphi^{(t)}$ comme un polyn\^ome en les variables $\wp(z)$, $\wp(Nz)$, $\frac{\wp'(Nz)}{\wp'(z)}$ de degr\'e partiels en $\wp(Nz)$ et $\frac{\wp'(Nz)}{\wp'(z)}$ respectivement major\'es par $3(L+t)\leq 6L$ et par $1$. On en d\'eduit l'existence d'une fraction rationnelle $G$, telle que $S_N(X)^{8L}G(X)$ soit un polyn\^ome \`a coefficients dans $\mathcal{O}_{\overline{K}}$ et v\'erifiant
\[ G(\wp(z))=\varphi^{(t)}(z)^2.\]
\noindent La fraction rationnelle $G(X)$ admet donc un z\'ero d'ordre sup\'erieur \`a $T_1$ aux points $X=\wp(u_1),\ldots,X=\wp(u_D)$. Notons $\Delta=\sum a_iX^i$ le polyn\^ome minimal unitaire sur $L$ de $\wp(u)$. Par hypoth\`ese sur $\wp(u)$, il est \`a coefficients entiers alg\'ebriques. Ainsi, il existe un polyn\^ome $H$ \`a coefficients $v$-entiers, tel que 
\begin{equation}\label{delta}
S_N(X)^{8L}G(X)=\Delta(X)^{T_1}H(X).
\end{equation}
\noindent Si $\pi_p$ est une uniformisante au-dessus de $p$ dans $K$ et $\pi_{p,F}$ une uniformisante de $p$ dans $F$, par le petit th\'eor\`eme de Fermat et le lemme \ref{nonramifie}, on obtient donc
\[\tau\left(\Delta\right)\left(\wp(\alpha_v u)\right)=\tau\left(\Delta\right)\left(\wp(u)^p\right)=\left(\Delta(\wp(u))\right)^p=0\ \mod \pi_{p,F}, \]
\noindent et ce, pour tout $\pi_{p,F}$ au-dessus de $\pi_p$. En substituant $\wp(\alpha_v u)$ \`a $X$ dans (\ref{delta}) et en appliquant $\tau$ aux coefficients des polyn\^omes $S$, $G$, $H$ et $\Delta$, on en d\'eduit que le membre de droite de cette \'egalit\'e transform\'ee par $\tau$ est d'ordre en $\pi_p$ sup\'erieur \`a $\frac{T_1}{e_p}$. Il reste maintenant \`a majorer l'ordre en $\pi_p$ de $\tau(S_N)(\wp(\alpha_v u)$. Or $S_N$ est \`a coefficients dans $K$, donc est, de m\^eme que $R_N$, invariant par $\tau$. De plus, 
\[\wp(N\alpha_v u)=\frac{R_N(\wp(\alpha_v u))}{S_N(\wp(\alpha_v u))}.\]
\noindent D'apr\`es le lemme \ref{premier}, les polyn\^omes $R_N$ et $S_N$ r\'eduits mod $\pi_p$ sont premiers entre eux. Autrement dit, l'un des deux nombres $R_N(\wp(\alpha_v u))$, $S_N(\wp(\alpha_v u))$ est une unit\'e de $\mathcal{O}_{\overline{K_v}}$. Si c'est $S_N(\wp(\alpha_v u))$, on a fini, sinon, c'est $R_N(\wp(\alpha_v u))$ et donc
\[\textnormal{ord}_{\pi_p}(S_N(\wp(\alpha_v u))=-\textnormal{ord}_{\pi_p}\wp(N\alpha_v u).\]
\noindent Cas 2 : Si $\wp(u)$ n'est pas un $v$-entier, on fait un changement de carte comme dans le b) du deuxi\`eme pas de Laurent \cite{laurent} : on effectue le changement de variable projectif 
\[t=-\frac{X}{Y}\ \ \textnormal{ et, }\ \ s=-\frac{1}{Y}.\]
\noindent Alors $s$ s'exprime en fonction du param\`etre local $t$ par une s\'erie enti\`ere $s(t)$ \`a coefficients $v$-entiers. On consid\`ere cette fois la fonction 
\[ \left(\wp'(z)\wp'(Nz)\right)^{-(L_1+t)}\varphi^{(t)}(z), \textnormal{ \`a la place de }\varphi^{(t)}(z)^2.\]
\noindent Elle s'\'ecrit comme un polyn\^ome en les variables $-\frac{\wp(z)}{\wp'(z)}$, $-\frac{1}{\wp(z)}$, $-\frac{\wp(Nz)}{\wp'(Nz)}$ et  $-\frac{1}{\wp'(Nz)}$. Il existe donc une s\'erie enti\`ere $G$, \`a coefficients $v$-entiers, telle que
\begin{equation}\label{g}
G(t)= \left(\wp'(z)\wp'(Nz)\right)^{-(L_1+t)}\varphi^{(t)}(z).
\end{equation}
\noindent Soient $\xi=-\frac{\wp(u)}{\wp'(u)}$ le param\`etre local associ\'e au point $P$ et $\Delta$ le polyn\^ome minimal unitaire de $\xi$ sur $F_v$. D'apr\`es l'hypoth\`ese, le nombre $\xi$ est dans l'id\'eal maximal $\overline{\mathfrak{m}_v}$, donc tous les coefficients de $\Delta$, sauf le coefficient dominant, sont divisibles par $\pi_p$. Le th\'eor\`eme de pr\'eparation de Weierstrass montre qu'il existe une s\'erie enti\`ere $H$, \`a coefficients $v$-entiers, telle que 
\[G(t)=\Delta(t)^{(T_1-t)}H(t).\]
\noindent Dans cette identit\'e on substitue $t=[\alpha_v](\xi)=\xi^p+\pi_p\psi(\xi)$ par le lemme \ref{prep}. On conclut alors comme dans le premier cas \Large{(}\normalsize il faut savoir majorer
\[\textnormal{ord}_{\pi_p}(\wp'(\alpha_v u)\wp'(N\alpha_vu))=\textnormal{ord}_{\pi_p} s([\alpha_v](\xi))+\textnormal{ord}_{\pi_p} s([N\alpha_v](\xi)).\]
\noindent Puisque $N$ est premier \`a $p$, ces derni\`eres quantit\'es sont toutes deux \'egales \`a
\[-\frac{3}{2}\textnormal{ord}_{\pi_p}\wp(N\alpha_v u) \textnormal{ (cf. paragraphes 3 et 6 de \cite{tate}).}\]
\noindent Rempla\c{c}ant dans (\ref{g}), $z$ par $\alpha_v u$ et $t$ par $[\alpha_v](\xi)$ on peut alors conclure\Large{)}\normalsize.\hfill$\Box$

\vspace{.3cm}

\noindent On passe maintenant \`a la proposition qui nous servira quand il y a beaucoup de grande ramification.

\vspace{.3cm}

\begin{prop}\label{r}Soient $L_2$ et $T_2$ deux entiers d'ordre de grandeur polynomial en $D$ et $\Lambda_2$ un entier strictement positif, tels que $L_2^2\geq DT_2\Lambda_2$. On consid\`ere la fonction $\varphi$ obtenue dans la proposition \ref{siegelth} avec $L=L_2$, $T=T_2$, $k=D\Lambda_2$ et avec 
\[\left\{u_1,\ldots,u_k\right\}:=\left\{\alpha_v u_i\ / i\in\left\{1,\ldots,D\right\}\textnormal{ et }v \textnormal{ d\'ecrivant un ensemble $\mathcal{P}_2^{\star}$ de cardinal }\Lambda_2\right\}.\]
\noindent Pour tout $t\leq L_2$, pour tout $v$ dans l'ensemble $\mathcal{P}_2^{\star}$ et pour tout $\tau\in\textnormal{Gal}(\overline{K}/K)$ tel que $\tau_{F}\in H_p$, $p/v\in\mathcal{P}_2^{\star}$, on a
\[\log \mid\tau(\varphi^{(t)})(\alpha_p u)\mid_v\leq -\frac{T_2}{2}\log p+8L_2\log\max\{1,\mid\wp(N\alpha_vu)\mid_v\}.\]
\end{prop}
\demo \noindent L\`a encore il y a deux cas selon que $\wp(u)$ est un $v$-entier ou non. On commence par le cas o\`u c'est un $v$-entier. On \'etend $v$ au corps $\overline{K}(X)$ en posant $\mid X\mid_v=1$ et on fait la m\^eme preuve que pr\'ec\'edemment, en montrant cette fois-ci que 
\[\tau\left(\Delta_{\alpha_v}\right)(\wp(\alpha_v u))=0\mod \pi_p,\]
\noindent o\`u $\Delta_{\alpha_v}$ est le polyn\^ome minimal de $\wp(\alpha_v u)$. Pour montrer ceci, on note $\Delta^{(p)}$ le polyn\^ome minimal de $\wp(u)$ o\`u l'on a \'elev\'e les coefficients \`a la puissance $p$. On a alors 
\begin{align*}
\tau\left(\Delta_{\alpha_v}\right)(\wp(\alpha_v u))	& = \tau\left(\Delta^{(p)}\right)(\wp(\alpha_v u))\mod \pi_p\textnormal{ par le petit th\'eor\`eme de Fermat,}\\
							& = \Delta^{(p)}(\wp(\alpha_v u))\mod\pi_p\textnormal{ par le lemme \ref{ramifie}},\\
							& =\Delta^{(p)}(\wp(u)^p)\mod\pi_p,\\
							& = \left(\Delta(\wp(u))\right)^p\mod\pi_p,\\
							& = 0.\\
\end{align*}
\noindent Si $\wp(u)$ n'est pas un $v$-entier, on utilise le m\^eme argument que dans le cas 2 de la proposition \ref{nr} pr\'ec\'edente pour conclure de la m\^eme fa\c{c}on.\hfill$\Box$

\section{Conclusion}\label{conclusion}

\subsection{Le cas du th\'eor\`eme \ref{thmoi}}

\noindent En notant $[\cdot]$ la partie enti\`ere, on pose $C$ une constante assez grande (de sorte que les in\'egalit\'es soient v\'erifi\'ees) ne d\'ependant que de $E/K$ et on pose 
\[N_1=\left[C^4\frac{(\log 2D)^6}{(\log\log 5D)^5}\right]\textnormal{ et }\ E=\left[C\left(\frac{\log 2D}{\log\log 5D}\right)^2\right].\]
\noindent Pour $p$ entre $N_1/2$ et $N_1$, le th\'eor\`eme de Chebotarev nous indique qu'il y a plus de $\Lambda=\left[\frac{C^4}{2}\left(\frac{\log 2D}{\log\log 5D}\right)^6\right]$ tels $p$. En notant $e_v$ l'indice de ramification de $v$ dans $F$, on a : soit il y a plus de $\Lambda_1=\Lambda/2$ nombres premiers $p$ ayant une place $v$ avec un $e_v\leq E$, soit il y a plus de $\Lambda_2=\Lambda/2$ nombres premiers $p$ ayant toutes les places $v$ avec un $e_v>E$. On va traiter chaque cas s\'epar\'ement et conclure dans chacun de ces deux cas.

\vspace{.3cm}

\noindent Cas 1 : il y a plein de $v$ ayant peu de ramification, \textit{i.e.}, il y a plus de $\Lambda_1=\Lambda/2$  nombres premiers $p$ ayant une place $v$ avec un $e_v\leq E$.

\vspace{.3cm}

\noindent Dans ce cas, on note $\mathcal{P}_1^{\star}$ le sous-ensemble de $\mathcal{P}^{\star}$ correspondant \`a $\Lambda_1$ et on introduit les param\`etres suivants :
\[L_1=\left[C^3D\frac{(\log 2D)^5}{(\log\log 5D)^6}\right],\ T_1=\left[C^{\frac{9}{2}}D\frac{(\log 2D)^7}{(\log\log 5D)^9}\right],\textnormal{ et, } T_1'=\left[C^3D\frac{(\log 2D)^4}{(\log\log 5D)^6}\right],\]
\noindent  et $N$ est un nombre premier tel que $\frac{1}{2}\sqrt{L_1}\leq N\leq \sqrt{L_1}$. 

\vspace{.3cm}

\begin{prop} Pour tout $p\in\mathcal{P}_1^{\star}$, pour tout $\tau$ \'etendant $\Phi_p^{-1}$ et pour tout $t\leq T_1'$, la fonction 
$\varphi^{(t)}$ de la proposition \ref{nr} s'annule en $\tau\wp(\alpha_v u)$.
\end{prop}
\demo En notant 
\[\zeta=\textnormal{N}_{F(P)/F}\left(\tau\left(\varphi^{(t)}\right)(\alpha_v u)\right),\]
on a gr\^ace \`a la proposition \ref{nr}, 
\[\log \mid\zeta\mid_v\leq -\frac{DT_1\log p}{2e_v}+8L_1\sum_{w/v}D_w\log\max\left(1,\mid\wp(N\alpha_vu\mid_w\right)...\]
\noindent On en d\'eduit
\[\log \mid\zeta\mid_v\leq -c_{12}\frac{DT_1\log p}{2E}+DL_1\left(c_{10}+N^2p\widehat{h}(P)\right).\]
\noindent Or par hypoth\`ese sur $\widehat{h}(P)$, on a $N^2p\widehat{h}(P)\leq N^2N_1\widehat{h}(P)\leq c_{11}$. En rempla\c{c}ant les param\`etres par leur valeur, on obtient donc 
\[\log \mid\zeta\mid_v\leq -C^{\frac{7}{2}}D^2\frac{(\log 2D)^{5}}{(\log\log 5D)^{6}}.\]
\noindent Ainsi, si $\zeta$ est non nul, on a 
\begin{equation}\label{min}
h(\zeta)=h(\zeta^{-1})\geq\frac{d_v}{d}\log\max\left\{1,\mid\zeta^{-1}\mid_v\right\}\geq C^{\frac{7}{2}}D^2\frac{(\log 2D)^{5}}{(\log\log 5D)^{6}}.
\end{equation}

\noindent Par ailleurs, un calcul classique (cf. par exemple \cite{davidhindry} p.50) permet d'\'ecrire
\[h(\zeta)\leq c_{15}DT_1'\log(T_1'+L_1)+c_{16}DL_1N^2p\widehat{h}(P)+c_{17}Dh_2(\varphi)\]
\noindent o\`u $h_2(\varphi)$ est donn\'ee par la proposition (\ref{siegelth}) :
\[h_2(\varphi)\leq  \frac{cDT_1}{(L_1+1)^2-DT_1}\left(L_1N^2\widehat{h}(P)+T_1\log(T_1+L_1)+T_1\log N+L_1\right)+\log L_1.\]
\noindent en rempla\c{c}ant les param\`etres par leur valeur, on obtient :
\[h(\zeta)\leq c_{15}C^3D^2\frac{(\log 2D)^{5}}{(\log\log 5D)^6}+c_{16}C^3D^2\frac{(\log 2D)^{5}}{(\log\log 5D)^6}+c_{17}C^3D^2\frac{(\log 2D)^{5}}{(\log\log 5D)^6},\]
\noindent Soit
\begin{equation}\label{maj}
h(\zeta)\leq c_{18}C^3D^2\frac{(\log 2D)^{5}}{(\log\log 5D)^6}.
\end{equation}
\noindent En comparant les in\'egalit\'es (\ref{min}) et (\ref{maj}), on obtient une contradiction pour $C$ suffisamment grand, ce qui conclut.\hfill$\Box$

\vspace{.3cm}

\noindent Puisque l'on travaille avec des $p\in \mathcal{P}^{\star}\subset\mathcal{P}$, on a $[F(\alpha_v (P)) :F]=D$. Donc, on obtient ainsi, en comptant les multiplicit\'es, au moins 
\begin{equation}\label{racine}
DT_1'\Lambda_1\geq \frac{1}{2}C^7D^2\frac{(\log 2D)^{10}}{(\log\log 5D)^{12}}
\end{equation}
\noindent racines. Or en utilisant la relation 
\[\wp(Nz)=\frac{R_N(\wp(z))}{S_N(\wp(z))},\]
\noindent on peut \'ecrire $\varphi$ sous la forme 
\[\varphi(z)=F(\wp(z)),\]
\noindent o\`u $F$ est une fraction rationnelle de degr\'e major\'e par 
\begin{equation}\label{fraction}
(N^2+1)L_1\leq 2L_1^2\leq 2C^6D^2\frac{(\log 2D)^{10}}{(\log\log 5D)^{12}}.
\end{equation}
\noindent En comparant (\ref{racine}) et (\ref{fraction}), on en d\'eduit que la fraction $F$ est identiquement nulle. Donc il en est de m\^eme pour $\varphi$, ce qui est absurde par le choix de $N$. Le th\'eor\`eme est donc d\'emontr\'e dans ce cas.

\vspace{.3cm}

\noindent Cas 2 : il y a plein de $v$ ayant beaucoup de ramification, \textit{i.e.}, il y a plus de $\Lambda_2=\Lambda/2$ nombres premiers $p$ ayant toutes les places $v$ avec un $e_v> E$.

\vspace{.3cm}

\noindent Dans ce cas, on note $\mathcal{P}_2^{\star}$ le sous-ensemble de $\mathcal{P}^{\star}$ correspondant \`a $\Lambda_2$ et on introduit les param\`etres suivants :
\[L_2=\left[C^{\frac{35}{8}}D\frac{(\log 2D)^{7}}{(\log\log 5D)^8}\right],\ T_2=\left[C^{\frac{9}{2}}D\frac{(\log 2D)^{7}}{(\log\log 5D)^9}\right]\textnormal{ et, } T_2'=\left[C^4D\frac{(\log 2D)^6}{(\log\log 5D)^8}\right],\]
\noindent et $N$ est un nombre premier tel que $\frac{1}{2}\sqrt{L_2}\leq N\leq \sqrt{L_2}$. 

\vspace{.3cm}

\begin{prop} Pour tout $p\in\mathcal{P}_2^{\star}$, pour tout $\tau$ tel que $\tau_F\in H_p$ et pour tout $t\leq T_2'$, la fonction 
$\varphi^{(t)}$ de la proposition \ref{r} s'annule en $\tau\wp(\alpha_v u)$.
\end{prop}
\demo  En notant 
\[\zeta=\textnormal{N}_{F(P)/F}\left(\tau\left(\varphi^{(t)}\right)(\alpha_v u)\right),\]
on a gr\^ace \`a la proposition \ref{r}, 
\[\log \mid\zeta\mid_v\leq -\frac{DT_2\log p}{2}+\frac{8L_2}\sum_{w/v}D_w\log\max\left(1,\mid\wp(N\alpha_vu\mid_w\right).\]
\noindent On en d\'eduit
\[\log \mid\zeta\mid_v\leq -c_{12}\frac{DT_2\log p}{2}+DL_2\left(c_{10}+N^2p\widehat{h}(P)\right).\]
\noindent Or par hypoth\`ese sur $\widehat{h}(P)$, on a $N^2p\widehat{h}(P)\leq N^2N_1\widehat{h}(P)\leq c_{11}$. En rempla\c{c}ant les param\`etres par leur valeur, on obtient donc 
\[\log \mid\zeta\mid_v\leq -c_{13}C^{\frac{9}{2}}D^2\frac{(\log 2D)^7}{(\log\log 5D)^8}.\]
\noindent Ainsi, si $\zeta$ est non nul, on a 
\begin{equation}\label{mi}
h(\zeta)=h(\zeta^{-1})\geq\frac{d_v}{d}\log\max\left\{1,\mid\zeta^{-1}\mid_v\right\}\geq c_{14}C^{\frac{9}{2}}D^2\frac{(\log 2D)^7}{(\log\log 5D)^8}.
\end{equation}
\noindent Par ailleurs, le m\^eme calcul que pr\'ec\'edemment permet d'\'ecrire
\[h(\zeta)\leq c_{15}DT_2'\log(T_2'+L_2)+c_{16}DL_2N^2p\widehat{h}(P)+c_{17}Dh_2(\varphi)\]
\noindent o\`u $h_2(\varphi)$ est donn\'ee par la proposition (\ref{siegelth}) :
\[h_2(\varphi)\leq  \frac{cD\Lambda_2T_2}{(L_2+1)^2-D\Lambda_2T_2}\left(L_2N^2\widehat{h}(P)+T_2\log(T_2+L_2)+T_2\log N+L_2\right)+\log L_2.\]
\noindent en rempla\c{c}ant les param\`etres par leur valeur, on obtient :
\[h(\zeta)\leq c_{15}C^4D^2\frac{(\log 2D)^7}{(\log\log 5D)^8}+c_{16}C^{\frac{35}{8}}D^2\frac{(\log 2D)^7}{(\log\log 5D)^8}+c_{17}C^{\frac{17}{4}}D^2\frac{(\log 2D)^7}{(\log\log 5D)^8},\]
\noindent soit,
\begin{equation}\label{ma}
h(\zeta)\leq c_{18}C^{\frac{17}{4}}D^2\frac{(\log 2D)^7}{(\log\log 5D)^8}.
\end{equation}
\noindent En comparant les in\'egalit\'es (\ref{mi}) et (\ref{ma}), on obtient une contradiction, ce qui conclut.\hfill$\Box$

\vspace{.3cm}

\noindent Puisque l'on travaille avec des $p\in \mathcal{P}^{\star}\subset\mathcal{P}$, on a $[F(\alpha_v (P)) :F]=D$. Donc, on obtient ainsi, en comptant les multiplicit\'es, au moins 
\begin{equation}\label{racin}
EDT_2'\Lambda_2\geq \frac{1}{2}C^9D^2\frac{(\log 2D)^{14}}{(\log\log 5D)^{16}}
\end{equation}
\noindent racines. Or en utilisant la relation 
\[\wp(Nz)=\frac{R_N(\wp(z))}{S_N(\wp(z))},\]
\noindent on peut \'ecrire $\varphi$ sous la forme 
\[\varphi(z)=F(\wp(z)),\]
\noindent o\`u $F$ est une fraction rationnelle de degr\'e major\'e par 
\begin{equation}\label{fractio}
(N^2+1)L_2\leq 2L_2^2\leq 2C^{\frac{35}{4}}D^2\frac{(\log 2D)^{14}}{(\log\log 5D)^{16}}.
\end{equation}
\noindent En comparant (\ref{racin}) et (\ref{fractio}), on en d\'eduit que la fraction $F$ est identiquement nulle. Donc il en est de m\^eme pour $\varphi$, ce qui est absurde par le choix de $N$. Le th\'eor\`eme est donc d\'emontr\'e dans ce cas. Il est donc d\'emontr\'e dans tous les cas.\hfill$\Box$

\subsection{Le cas du th\'eor\`eme \ref{th2}}

\noindent Pour prouver le th\'eor\`eme \ref{th2}, on fait essentiellement la m\^eme preuve que pour le th\'eor\`eme \ref{thmoi} : on peut faire les m\^emes r\'eductions et on a uniquement besoin de la partie non-ramifi\'ee de la preuve pr\'ec\'edente. La seule chose qui change est le choix des param\`etres permettant de conclure.

\vspace{.3cm}

\noindent En notant $[\cdot]$ la partie enti\`ere et $C$ une constante assez grande (de sorte que les in\'egalit\'es soient v\'erifi\'ees) ne d\'ependant que de $E/K$ et de $c_0$, on pose $E=1$ et 
\[N_1=\left[3C^2\frac{(\log 2D)^2}{\log\log 5D}\right].\]
\noindent Pour $p$ entre $N_1/2$ et $N_1$, le th\'eor\`eme de Chebotarev nous indique qu'il y a plus de $\Lambda=\left[\frac{C^2}{2}\left(\frac{\log 2D}{\log\log 5D}\right)^2\right]$ tels $p$. On note $\mathcal{P}_1^{\star}$ le sous-ensemble de $\mathcal{P}$ correspondant \`a $\Lambda$ et on introduit les param\`etres suivants :
\[L_1=\left[C^2D\frac{\log 2D}{\log\log 5D}\right],\ T_1=\left[2CD\frac{\log 2D}{\log\log 5D}\right],\textnormal{ et, } T_1'=\left[\frac{C^2}{2}D\right],\]
\noindent et $N$ est un nombre premier tel que $\frac{1}{2}\sqrt{L_2}\leq N\leq \sqrt{L_2}$. 

\vspace{.3cm}

\noindent Le m\^eme argument qu'au cas 1 du paragraphe pr\'ec\'edent nous permet alors de conclure.

\vspace{.3cm}

\begin{rem} \textnormal{Notons que bien que la preuve de ce th\'eor\`eme \ref{th2} soit moralement la m\^eme que celle du th\'eor\`eme de \cite{laurent}, le fait d'utiliser un lemme de Siegel absolu conduit \`a un choix diff\'erent des param\`etres pour faire fonctionner l'\'etape d'extrapolation.}
\end{rem}

\section{Application du th\'eor\`eme \ref{thmoi}}

\noindent En fait une version affaiblie du th\'eor\`eme \ref{thmoi} suffit d\'ej\`a. Pr\'ecis\'ement, on utilisera le corollaire suivant :

\begin{cor}\label{month}Soient $E/K$ une courbe elliptique \`a multiplication complexe et $\varepsilon$ un r\'eel strictement positif. On note $K^{\textnormal{ab}}$ la cl\^oture ab\'elienne de $K$. Il existe une constante strictement positive $c(E/K,\varepsilon)$ telle que 
\[\forall P\in E(\overline{K})\backslash E_{\textnormal{tors}},\ \ \widehat{h}(P)\geq\frac{c(E/K,\varepsilon)}{D^{1+\varepsilon}},\]
\noindent o\`u $D=[K^{\textnormal{ab}}(P):K^{\textnormal{ab}}]$.
\end{cor}

\vspace{.3cm}

\noindent On reprend les notations de l'article \cite{viada}, et on va montrer comment \'eviter la seconde partie de la preuve (parties 4.3 et 4.4 et 6 de l'article \cite{viada}).

\vspace{.3cm}

\noindent Soient $n\geq r\geq 0$ et $P\in S_{n-r}(C)$. On note $K(P)$ le corps de d\'efinition de $P$. On note $x_i : C\rightarrow E$ les applications coordonn\'ees d\'efinies par la composition de l'immersion ferm\'ee $C\hookrightarrow E^n$ et de la $i$-\`eme projection $E^n\rightarrow E$. Comme $E$ est \`a multiplication complexe, on sait que son anneau d'endomorphismes est un ordre $\mathcal{O}=\mathbb{Z}+\tau\mathbb{Z}$ dans un corps quadratique imaginaire. On d\'efinit alors $n$ morphismes suppl\'ementaires : $\forall 1\leq i\leq n,$ $x_{n+i}=\tau x_i$. On note 
\[\Gamma=\left\langle x_1,\ldots, x_n\right\rangle_{\textnormal{End}(E)}\]
le \textit{coordinate module} (d\'efinition de \cite{viada} p.51) : c'est le $\mathbb{Z}$-module engendr\'e par les $x_i$ avec $1\leq i\leq 2n$. On consid\`ere par ailleurs le $\mathbb{Z}$-module (qui est de rang $2r$ comme il suit du lemma 2. de \cite{viada}) 
\[\Gamma_P:=\left\langle x_1(P),\ldots,x_{2n}(P)\right\rangle_{\textnormal{End}(E)}.\]
\noindent Dans sa proposition 2. de \cite{viada}, Viada montre que $K\left((\Gamma_P)_{\textnormal{tors}}\right)\subset K(P)$ et elle montre \'egalement qu'il existe des \'el\'ements $\mathbb{Z}$-lin\'eairement ind\'e\-pendants $g_1,\ldots,g_{2r}$ de $\Gamma_P$, d\'efinis sur $K(P)$ et engendrant la partie libre de $\Gamma_P$. Ainsi on peut \'ecrire pour tout $1\leq i\leq 2n$,
\[x_i(P)=\sum_{j=1}^{2r}a_{ij}g_j +T_i\]
\noindent o\`u $T_i$ est un point de torsion. On pose $\nu_j=(a_{1j},\ldots,a_{nj})$ et on pose $\mid\nu_j\mid=\max_i\mid a_{ij}\mid$. Avec ces notations on a l'in\'egalit\'e (19) de \cite{viada} :
\begin{equation}\label{19}
\prod_{i=1}^{2r}\widehat{h}(g_i)\ll \prod_{i=1}^{2r}\mid\nu_i\mid^{-2}.
\end{equation}
\noindent Dans son corollary 1. Viada obtient alors l'in\'egalit\'e
\begin{equation}\label{23}
d\ll \left(NR\prod_{i=1}^{2r}\mid\nu_i\mid\right)^{\frac{1}{n-r}}
\end{equation}
\noindent o\`u $d=[K(P):K]$ et $N$ et $R$ sont deux entiers tels que $\left(\Gamma_P\right)_{\textnormal{tors}}\simeq \mathbb{Z}/N\mathbb{Z}\times\mathbb{Z}/R\mathbb{Z}$. (Ces entiers existent par la proposition 2. de \cite{viada}.)

\vspace{.3cm}

\noindent Dans son corollary 2. Viada obtient enfin les in\'egalit\'es 
\[ (NR)^{1-\varepsilon}\ll d \ll (NR)^{\frac{1}{n-r-1-\varepsilon}}.\]
\noindent Ainsi $d$ est born\'e en fonction de $N$ et si $n-r\geq 3$ on obtient l'in\'egalit\'e 
\[(NR)^{1-\varepsilon}\ll (NR)^{\frac{1}{2-\varepsilon}}\]
\noindent Ce qui permet de borner $N$ et donc de conclure. Dans le cas o\`u $n-r=2$, autrement dit le cas qui nous int\'eresse r\'eellement, on obtient juste 
\[(NR)^{1-\varepsilon}\ll (NR)^{\frac{1}{1-\varepsilon}}\]
\noindent ce qui ne permet malheureusement pas de conclure, d'o\`u la n\'ecessit\'e d'une seconde \'etape assez technique dans l'article \cite{viada}. On montre maintenant, et c'est l\`a la nouveaut\'e, comment conclure dans le cas g\'en\'eral en utilisant notre corollaire \ref{month} et en r\'eutilisant ce qui a \'et\'e fait jusqu'\`a pr\'esent. On se place d\'esormais dans le cas o\`u $n-r=2$. On note $K_N=K\left(\left(\Gamma_P\right)_{\textnormal{tors}}\right)$. On peut toujours supposer que $K=K(j(E))$, donc $K_N/K$ est une sous-extension ab\'elienne de l'extension ab\'elienne $K\left(E[N]\right)/K.$ On pose $D=\left[K(P):K_N\right]$. On a, en utilisant toujours le corollary 2. de \cite{viada},
\begin{equation}\label{*}
D=\frac{d}{[K_N:K]}\ll  (NR)^{\frac{1}{1-\varepsilon}}(NR)^{-(1-\varepsilon)}\leq (NR)^{3\varepsilon}
\end{equation}
\noindent si $\varepsilon$ est suffisament petit. Par ailleurs, les points $g_1$,\ldots, $g_{2r}$ sont des points d'ordre infini de $E(K(P))$. En appliquant le corollaire \ref{month} puis l'in\'egalit\'e (\ref{*}), on obtient
\begin{equation}\label{*1}
\prod_{i=1}^{2r}\widehat{h}(g_i)\gg D^{-2r-2r\varepsilon}\gg (NR)^{-6r\varepsilon(1+\varepsilon)}\gg (NR)^{-12n\varepsilon}.
\end{equation}
 
\noindent On a ainsi
\begin{align*}
(NR)^{1-\varepsilon}\ll d	& \ll \left(NR\prod_{i=1}^{2r}\mid\nu_i\mid\right)^{\frac{1}{2}}\text{ par l'in\'egalit\'e (\ref{23})}\\
			& \ll	(NR)^{\frac{1}{2}}\prod_{i=1}^{2r}\widehat{h}(g_i)^{-\frac{1}{4}} \text{ par l'in\'egalit\'e (\ref{19})}\\
			& \ll (NR)^{\frac{1}{2}+3n\varepsilon} \text{ par l'in\'egalit\'e (\ref{*1})}
\end{align*}
\noindent Ceci permet de conclure la preuve du th\'eor\`eme en prenant $\varepsilon$ assez petit.

\vspace{1cm}

\noindent \textbf{Adress :} RATAZZI Nicolas

Universit\'e Paris 6
Institut de Math\'ematiques

Projet Th\'eorie des nombres 

Case 247 

4, place Jussieu 

75252 Paris Cedex 05

FRANCE 

email : ratazzi@math.jussieu.fr

\end{document}